\documentclass{gtart_h}

\def\ifplaintex{\expandafter\ifx\csname documentclass\endcsname\relax}

\def\gtp{{\mathsurround=0pt\it $\cal G\mskip-2mu$eometry \&\ 
$\cal T\!\!$opology $\cal P\!$ublications}}  

\def\recd{{\small Received:\qua\receiveddate\ifx\reviseddate\relax
\else\qquad Revised:\qua\reviseddate\fi\par}} 

\def\lognumber#1{\def\thelognumber{#1}}
\def\volumenumber#1{\def\thevolumenumber{#1}}
\def\volumeyear#1{\def\thevolumeyear{#1}}
\def\papernumber#1{\def\thepapernumber{#1}}
\def\pagenumbers#1#2{\def\startpage{#1}\def\finishpage{#2}}
\def\published#1{\def\publishdate{#1}}

\def\received#1{\def\receiveddate{#1}}

\def\accepted#1{\def\accepteddate{#1}}

\def\asciiauthors#1{\def\theasciiauthors{#1}}
\def\asciiaddress#1{\def\theasciiaddress{#1}}
\def\asciiemail#1{\def\theasciiemail{#1}}

\def\coverauthors#1{\def\thecoverauthors{#1}}
\long\def\asciiabstract#1{\long\def\theasciiabstract{#1}}

\let\\\par\let\thelognumber\relax\let\thevolumenumber\relax
\let\thepapernumber\relax\let\thevolumeyear\relax\let\startpage\relax
\let\finishpage\relax\let\publishdate\relax\let\receiveddate\relax
\let\reviseddate\relax\let\accepteddate\relax\let\theasciititle\relax
\let\theasciiauthors\relax\let\theasciiaddress\relax
\let\theasciiabstract\relax

\let\thecoverauthors\relax\let\theasciiemail\relax

\ifplaintex
\font\logobig=cmssbx10 scaled 3836
\font\logomed=cmssbx10 scaled 2557
\else
\font\logobig=cmssbx10 scaled 4200
\font\logomed=cmssbx10 scaled 2800
\fi

\long\def\makeagttitle{   
\count0=\startpage
\agt\hfill      
\hbox to 45truept{\vbox to 0pt{\vglue -13truept{\logomed A\kern -.37em{\logobig 
T}\kern -.38em G}\vss}\hss}
\break
{\small Volume \thevolumenumber\ (\thevolumeyear)
\startpage--\finishpage\nl
Published: \publishdate}

\vglue .25truein

{\parskip=0pt\leftskip 0pt plus
1fil\def\\{\par\smallskip}{\Large\bf\thetitle}\par\medskip} \vglue
0.05truein

{\parskip=0pt\leftskip 0pt plus 1fil\def\\{\par}{\sc\theauthors}
\par\medskip}%
 
\vglue 0.03truein 

{\small\leftskip 25truept\rightskip 25truept{\bf Abstract}\stdspace\theabstract

{\bf AMS Classification}\stdspace\theprimaryclass
\ifx\thesecondaryclass\relax\else; \thesecondaryclass\fi\par
{\bf Keywords}\stdspace \thekeywords\par}\vglue 7truept

}   

\ifplaintex
\font\phead=cmsl9 scaled 950
\font\pnum=cmbx10 scaled 913
\font\pfoot=cmsl9 scaled 950
\headline{\vbox to 0pt{\vskip -4.5mm\line{\small\phead\ifnum
\count0=\startpage ISSN 1472-2739 (on-line) 1472-2747 (printed)
\hfill {\pnum\folio}\else\ifodd\count0\def\\{ }%
\ifx\theshorttitle\relax\thetitle\else\theshorttitle\fi\hfill{\pnum\folio}
\else\def\\{ and }{\pnum\folio}\hfill\ifx\theshortauthors\relax\theauthors
\else\theshortauthors\fi\fi\fi}\vss}}
\footline{\vbox to 0pt{\vglue 0mm\line{\small\pfoot\ifnum\count0=\startpage
\copyright\ \gtp\hfill\else
\agt, Volume \thevolumenumber\ (\thevolumeyear)\hfill\fi}\vss}}
\else
\font=cmsl9 scaled 1050
\font\lnum=cmbx10 
\font=cmsl9 scaled 1050
\makeatletter
\def\@oddhead{{\small\lhead\ifnum\count0=\startpage ISSN 1472-2739 
(on-line) 1472-2747 (printed)\hfill {\lnum\number\count0}\else\ifodd\count0
\def\\{ }\ifx\theshorttitle\relax \thetitle \else\theshorttitle\fi\hfill
{\lnum\number\count0}\else\def\\{ and }{\lnum\number\count0}
\hfill\ifx\theshortauthors\relax 
\theauthors\else\theshortauthors\fi\fi\fi}}\def\@evenhead{\@oddhead}
\def\@oddfoot{\small\lfoot\ifnum\count0=\startpage\copyright\ \gtp\hfill\else
\agt, Volume \thevolumenumber\ (\thevolumeyear)\hfill\fi}
\def\@evenfoot{\@oddfoot}
\makeatother
\fi
\let\maketitlepage\makeagttitle
\let\makeshorttitle\maketitlepage
\let\maketitle\maketitlepage

\newwrite\gtoutfile
\long\gdef\makeheadfile{  
{\def\\{, }\def\s{ }
\immediate\openout\gtoutfile head.xxx
\immediate\write\gtoutfile{Proxy-for: \ifx\theasciiauthors\relax
\theauthors\else\theasciiauthors\fi\s<\ifx\theasciiemail\relax\theemail\else\theasciiemail\fi>}
\immediate\write\gtoutfile{\noexpand\\}
\immediate\write\gtoutfile{Authors: \ifx\theasciiauthors\relax
\theauthors\else\theasciiauthors\fi}
{\def\\{ }\immediate\write\gtoutfile{Title: \ifx\theasciititle\relax
\thetitle\else\theasciititle\fi}}
\immediate\write\gtoutfile{Subj-class: GT or SG, GR etc}
\immediate\write\gtoutfile{MSC-class: \theprimaryclass\ifx\thesecondaryclass\relax\else, \thesecondaryclass\fi}
\immediate\write\gtoutfile{Journal-ref: Algebr. Geom. Topol. \thevolumenumber\s
(\thevolumeyear) \startpage-\finishpage}
\immediate\write\gtoutfile{Comments: Published by Algebraic and
Geometric Topology at}
\immediate\write\gtoutfile{\s\s\s  http://www.maths.warwick.ac.uk/agt/AGTVol\thevolumenumber/agt-\thevolumenumber-\thepapernumber.abs.html}
\immediate\write\gtoutfile{\noexpand\\}
\immediate\write\gtoutfile{}
\ifx\theasciiabstract\relax
\immediate\write\gtoutfile{\theabstract}\else
\immediate\write\gtoutfile{\theasciiabstract}\fi
\immediate\write\gtoutfile{}
\immediate\write\gtoutfile{\noexpand\\}
\immediate\write\gtoutfile{}
\immediate\closeout\gtoutfile}}  

\def\maketitlepage{\makeagttitle\makeheadfile}
\let\makeshorttitle\maketitlepage
\let\maketitle\maketitlepage

\lognumber{68}
\volumenumber{5}
\volumeyear{2005}
\papernumber{68}
\pagenumbers{1677}{1710}
\received{13 June 2005} 
\accepted{28 November 2005}
\published{7 December 2005}

\usepackage{amsmath,amssymb,psfrag,graphicx,subfigure}
\usepackage[all]{xy}

\def\psfraga <#1,#2> #3#4{%
\psfrag {#3}{\smash{\rlap{\kern #1 \raise #2\hbox{#4}}}}}

\def\fref#1{\hyperlink{#1anchor}{\ref*{#1}}}
\def\figref#1{\hyperlink{#1anchor}{Figure~\ref*{#1}}}
\def\tabref#1{\hyperlink{#1anchor}{Table~\ref*{#1}}}
\def\anchor#1{\noindent\hypertarget{#1anchor}{\smash{$\phantom{99}$}}}

\newtheorem{thm}{Theorem}[section]
\newtheorem{cor}[thm]{Corollary}
\newtheorem{lem}[thm]{Lemma}
\newtheorem{prop}[thm]{Proposition}
\theoremstyle{definition}
\newtheorem*{rems}{Remarks}
\newtheorem{rem}[thm]{Remark}

\newcommand\leftidx[3]{%
  {\vphantom{#2}}#1#2#3%
}

\newcommand{\id}{\mathrm{id}}

\newcommand{\diag}{\mathcal{D}iag}
\newcommand{\PB}{\mathrm{PB}}
\newcommand{\RPB}{\mathrm{RPB}}
\newcommand{\qdiagS}{\overline{\mathcal{D}iag}\,\!^S}
\newcommand{\qddS}{\overline{\mathcal{D}}\,\!^S}
\newcommand{\rhand}{\mathrm{RHand}}
\newcommand{\rstl}{\mathrm{RSL}}
\newcommand{\conv}{\mathrm{Conv}}

\newcommand{\cc}{\mathcal{C}}

\newcommand{\dd}{\mathcal{D}}

\newcommand{\trait}{\nobreakdash-\hspace{0pt}}
\newcommand{\N}{\mathbb{N}}

\newcommand{\un}{\mathbf{1}}

\newcommand{\Ob}{\mathrm{Ob}}

\newcommand{\End}{\mathrm{End}}
\newcommand{\Hom}{\mathrm{Hom}}

\newcommand{\ev}{\mathrm{ev}}

\newcommand{\coev}{\mathrm{coev}}

\newcommand{\dual}[1]{\leftidx{^\vee}{\!#1}{}}

\newcommand{\scaledraw}[1]{\includegraphics[scale=.8]{#1.eps}}
\newcommand{\scaleraisedraw}[2]{\raisebox{-#1\height}{\includegraphics[scale=.8]{#2.eps}}}
\newcommand{\pdraw}[1]{\raisebox{-.2\height}{\includegraphics[scale=.9]{#1.eps}}}

\newcommand{\TabTable}{\rule{0pt}{2.6ex}\rule[-1.2ex]{0pt}{0pt}}

\begin{document}

\title{Hopf diagrams and quantum invariants}
\authors{Alain Brugui\`eres\\Alexis Virelizier}
\coverauthors{Alain Brugui\noexpand\`eres\\Alexis Virelizier}
\asciiauthors{Alain Bruguieres and Alexis Virelizier}
\gtemail{\mailto{bruguier@math.univ-montp2.fr}{\rm\qua
and\qua}\mailto{virelizi@math.berkeley.edu}}
\asciiemail{bruguier@math.univ-montp2.fr, virelizi@math.berkeley.edu}
\address{I3M, Universit\'e Montpellier II, 34095 Montpellier Cedex 5, France
\\{\rm and}\\
Department of Mathematics, University of California, Berkeley CA 94720,
USA}

\asciiaddress{I3M, Universite Montpellier II, 34095 Montpellier Cedex 5, France
\\and\\
Department of Mathematics, University of California, Berkeley CA 94720,
USA}

\begin{abstract}
The Reshetikhin-Turaev invariant, Turaev's TQFT, and many related
constructions rely on the encoding of certain tangles ($n$-string
links, or ribbon $n$-handles) as $n$-forms on the coend of a ribbon
category. We introduce the monoidal category of Hopf diagrams, and
describe a universal encoding of ribbon string links as Hopf
diagrams. This universal encoding is an injective monoidal functor and
admits a straightforward monoidal retraction. Any Hopf diagram with
$n$ legs yields a $n$-form on the coend of a ribbon category in a
completely explicit way. Thus computing a quantum invariant of a
$3$-manifold reduces to the purely formal computation of the
associated Hopf diagram, followed by the evaluation of this diagram in
a given category (using in particular the so-called Kirby elements).
\end{abstract}

\asciiabstract{%
The Reshetikhin-Turaev invariant, Turaev's TQFT, and many related
constructions rely on the encoding of certain tangles (n-string links,
or ribbon n-handles) as n-forms on the coend of a ribbon category. We
introduce the monoidal category of Hopf diagrams, and describe a
universal encoding of ribbon string links as Hopf diagrams. This
universal encoding is an injective monoidal functor and admits a
straightforward monoidal retraction. Any Hopf diagram with n legs
yields a n-form on the coend of a ribbon category in a completely
explicit way. Thus computing a quantum invariant of a 3-manifold
reduces to the purely formal computation of the associated Hopf
diagram, followed by the evaluation of this diagram in a given
category (using in particular the so-called Kirby elements).}


\maketitle

\section*{Introduction}
\addcontentsline{toc}{section}{Introduction}

In 1991, Reshetikhin and Turaev \cite{RT2} introduced  a new $3$-manifold invariant. The construction proceeds in two steps:
representing a $3$-manifold by surgery along a link and then {\it coloring} the link to obtain a scalar invariant. Here,
{\it colors} are (linear combinations of) simple representations of a quantum group at a root of unity. Since then, this
construction has been re-visited many times.

In particular, Turaev \cite{Tur2} introduced the notion of a modular category, which is (after innocuous additivisation and
karoubianisation, see \cite{Brug2}) an abelian semisimple ribbon category satisfying a finiteness and a non-degeneracy
condition, and showed that such a category defines a $3$-manifold invariant, and indeed a TQFT. In this approach, colors are
simple objects of the category.

Following these ideas, a more general approach on quantum invariants of $3$\trait manifolds has been subsequently developed,
see \cite{Lyu2} and more recently \cite{LyuKer,Vir}. It avoids in particular the semisimplicity condition. Let us briefly
outline it: the initial data used to construct the invariants is a ribbon category $\cc$ endowed with a coend $A=\int^{X \in
\cc} \dual{X} \otimes X$. Let $L$ be a framed $n$-link. We can always present $L$ as the closure of some ribbon $n$-string
link $T$. By using the universal property of the coend $A$, to such a string link $T$ is associated a $n$-form on $A$, that
is, a morphism $T_\cc \co A^{\otimes n} \to \un$ in $\cc$. Given a morphism $\alpha\co \un \to A$ (which plays here the role
of the color), set:
\begin{equation*}
\tau_\cc(L;\alpha)=T_\cc \circ \alpha^{\otimes n} \in \End_\cc(\un).
\end{equation*}
A \emph{Kirby element} of $\cc$, as defined in \cite{Vir}, is a morphism $\alpha\co \un \to A$ such that, for all framed
link~$L$, $\tau_\cc(L;\alpha)$ is well-defined and invariant under isotopies of $L$ and under 2-handle slides. In this case,
by Kirby's theorem \cite{Ki} and under some invertibility condition, the invariant $\tau_\cc(L;\alpha)$ can be normalized to
an invariant of 3-manifolds.

At this stage, two main questions naturally arise. Firstly, how to recognize the Kirby elements of a ribbon category? And
secondly, how to compute the forms~$T_\cc$ obtained via the universal property?

Concerning the first question, recall that the coend $A$ has a structure of a Hopf algebra in the category $\cc$, see
\cite{Maj1,Lyu1}. This means in particular that $A$ is endowed with a product $\mu_A \co A\otimes A \to A$, a coproduct
$\Delta_A \co  A \to A \otimes A$, and an antipode $S_A \co  A \to A$ which satisfy the same axioms as those of a Hopf
algebra except one has to replace the usual flip map with the braiding of $\cc$. If $A$ admits a two-sided integral
$\lambda\co \un \to A$, then $\lambda$ is a Kirby element and the corresponding $3$\trait manifold invariant is that of
Lyubashenko~\cite{Lyu2}. More generally, if a morphism $\alpha\co \un \to A$ in $\cc$ satisfies:
\begin{equation*}
(\id_A \otimes \mu_A)(\Delta_A \otimes \id_A)(\alpha \otimes \alpha)=\alpha \otimes \alpha \quad \text{and} \quad S_A
\alpha=\alpha,
\end{equation*}
then it is a Kirby element, see \cite{Vir}. In particular, any Reshetikhin-Turaev invariant computed from a semisimple
sub-quotient of $\cc$ can be defined directly via a Kirby element of $\cc$ satisfying this last equation.

\begin{figure}[ht!]\anchor{figintro}
   \begin{center}
$D=$\,
\begin{minipage}[c]{1.52cm}\scalebox{.8}{\psfrag{S}[c][c]{$S$}\psfrag{D}[c][c]{$\Delta$}\psfrag{w}[c][c]{$\omega$}
\includegraphics{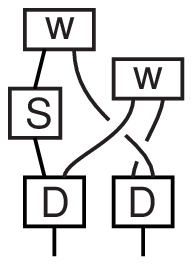}}\end{minipage}
$\rightsquigarrow$\,\;
\begin{minipage}[c]{2.56cm}\scalebox{.8}{\psfrag{S}[c][c]{$S$}\psfrag{D}[c][c]{$\Delta$}\psfrag{w}[c][c]{$\omega$}
\includegraphics{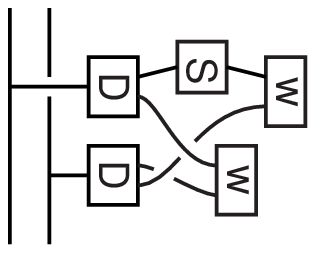}}\end{minipage}
$\rightsquigarrow$ $\;\psi(D)=$\,
\begin{minipage}[c]{4.3cm}\scalebox{.8}{\psfrag{=}[c][c]{$\sim$}\includegraphics{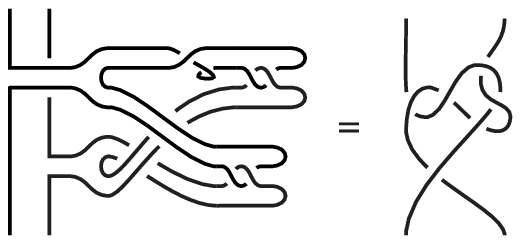}}\end{minipage}
   \end{center}
     \caption{From Hopf diagrams to ribbon string links}
     \label{figintro}
\end{figure}

The main motivation of the present paper is to answer the second question. Given a ribbon string link $T$, we express
$T_\cc$ in terms of some structural morphisms of $A$ (avoiding the product). This can be done by means of a universal
formula, that is, independently of $\cc$. To this end, we introduce the notion of \emph{Hopf diagrams}. They are braided
planar diagrams (with inputs but no output) obtained by stacking boxes such as $\Delta\!=\!\pdraw{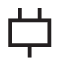}\,$,
$S\!=\!\pdraw{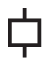}\,$, and $\omega\!=\!\pdraw{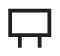}\,$. These are submitted to relations corresponding to those of a
coproduct, an antipode, and a Hopf pairing. To each Hopf diagram~$D$ with $n$ inputs is associated a ribbon $n$-string link
$\psi(D)$ as in \figref{figintro}.
We construct a category $\qdiagS$, whose objects are non-negative integers and morphisms are Hopf diagrams up to a certain
equivalence relation, in such a way that $\psi$ factors as a monoidal functor $\overline{\psi}$ from $\qdiagS$ to the
category $\rstl$ of ribbon string links. Moreover we construct a monoidal functor $\Psi\co\rstl \to\qdiagS$ which is right
inverse to $\overline{\psi}$, that is, such that $\overline{\psi}\circ \Psi= 1_{\rstl}$. This is the main result of this
paper (Theorem \ref{rstl2diag}). This monoidal functor may be viewed as a formal encoding of ribbon string links. More
precisely, given a ribbon category $\cc$ endowed with a coend $A=\int^{X \in \cc} \dual{X} \otimes X$, there is a canonical
map $E_\cc\co \{\text{Hopf diagrams with $n$ inputs} \} \to \Hom_\cc(A^{\otimes n}, \un)$ such that $T_\cc= E_\cc \circ
\Psi(T)$ for any ribbon $n$-string link  $T$ (see Theorem \ref{Evaldiag}).

In Section \ref{sect-Hopf}, we define the the monoidal category $\diag^S$ of Hopf diagrams, as a convolution category. The
category of Hopf diagrams comes in two versions: with, or without antipode. Both versions are isomorphic, however. In
Section~\ref{sect-hand}, we review the monoidal category $\rstl$ of ribbon string links and the related monoidal category of
ribbon handles. These categories are isomorphic. In Section~\ref{sect-hopf2hand}, we construct a monoidal functor
$\psi\co\diag^S \to \rstl$. In Section~\ref{sect-hand2hopf}, we define a category $\qdiagS$, as a quotient of $\diag^S$ by
certain new relations, in such a way that $\psi$ induces a monoidal functor $\overline{\psi}\co \qdiagS \to \rstl$ which
admits a right inverse. Finally, in Section \ref{sect-qinv}, given a ribbon category~$\cc$ endowed with a coend $A$, we
explain how to represent the category of Hopf diagrams into $\cc$ by using some structural morphisms of $A$. Moreover, we
give a general criterion, using Hopf diagrams,
for recognizing Kirby elements.\\

Unless otherwise specified, monoidal categories are assumed to be strict.\\

\rk{Acknowledgements} This work has been partially funded by the CNRS-NSF project n$^\circ$17149 `Algebraic and
homologic methods in low dimensional topology'. The second author thanks the Max-Planck-Institut f\"{u}r Mathematik  in Bonn,
where part of this work was carried out, for its support and hospitality.

\section{Hopf diagrams}\label{sect-Hopf}

In this section, we construct the categories $\diag$ and $\diag^S$ (as convolution categories) which are shown to be
isomorphic. They are preliminary versions of the category of Hopf diagrams.

\subsection{Categorical (co)algebras}
Recall that an \emph{algebra} in a monoidal category is an object $A$ endowed with morphisms $\mu\co A \otimes A \to A$ (the
\emph{product}) and  $\eta\co  A \to \un$ (the \emph{unit}) which satisfy:
\begin{equation*}
\mu(\id_A \otimes \mu)=\mu(\mu \otimes \id_A) \text{ \quad and \quad} \mu(\eta \otimes \id_A)=\id_A=\mu(\id_A \otimes \eta).
\end{equation*}
Dually, a \emph{coalgebra} in a monoidal category is an object $C$ endowed with morphisms $\Delta\co C \to C \otimes C$ (the
\emph{coproduct}) and $\varepsilon\co  C \to \un$ (the \emph{counit}) which satisfy:
\begin{equation*}
(\id_C \otimes \Delta)\Delta=(\Delta \otimes \id_C)\Delta \text{ \quad and \quad} (\varepsilon \otimes
\id_C)\Delta=\id_C=(\id_C \otimes \varepsilon)\Delta.
\end{equation*}

Note that the unit object $\un$ of a monoidal category $\cc$ is both an algebra (with $\mu=\id_\un=\eta$) and a coalgebra
(with $\Delta=\id_\un=\varepsilon$) in $\cc$ (recall that $\un=\un \otimes \un$).

When the category $\cc$ is braided with braiding $\tau$, the monoidal product $A\otimes A'$ of two algebras $A$ and $A'$ in
$\cc$ is an algebra in $\cc$ with unit $\eta \otimes \eta'$ and product $(\mu \otimes \mu')(\id_A \otimes \tau_{A,A'}
\otimes \id_{A'})$. In particular, for any non-negative integer $n$, $A^{\otimes n}$ is an algebra in~$\cc$. Likewise, the
monoidal product $C\otimes C'$ of two coalgebras $C$ and $C'$ in a braided category $\cc$ is a coalgebra in $\cc$ with
counit $\varepsilon \otimes \varepsilon'$ and coproduct $(\id_C \otimes \tau_{C,C'} \otimes \id_{C'})(\Delta \otimes
\Delta')$. In particular, for any non-negative integer $n$, $C^{\otimes n}$ is a coalgebra in~$\cc$.

\subsection{The convolution product}\label{convoprodsect}

Let $\cc$ be a monoidal category, $(A,\mu,\eta)$ be an algebra in $\cc$, and $(C,\Delta,\varepsilon)$ be a coalgebra in
$\cc$. The \emph{convolution product} of two morphisms $f,g \in \Hom_\cc(C,A)$ is the morphism $f\star g=\mu(f \otimes
g)\Delta \in \Hom_\cc(C,A) $. This makes the set $\Hom_\cc(C,A)$ a monoid with unit $\eta\varepsilon\co C \to A$.

\subsection{Convolution categories}\label{convcat}
Let $\cc$ be a braided category, $A$ be an algebra in~$\cc$, and $C$ be a coalgebra in~$\cc$. Let us define the
\emph{convolution category} $\conv_\cc(C,A)$ as follows: the objects of $\conv_\cc(C,A)$ are the non-negative integers $\N$.
For $m,n \in \N$, the set of morphisms from $m$ to $n$ is empty if $m \neq n$ and is the monoid $\Hom_\cc(C^{\otimes n},A)$
endowed with the convolution product if $m=n$ (recall indeed that $C^{\otimes n}$ is a coalgebra in $\cc$). In particular,
the identity of an object $n\in \N$ is:
\begin{equation*}
\id_n=\eta \varepsilon^{\otimes n}\co C^{\otimes n} \to A,
\end{equation*}
and the composition of two endomorphisms $f,g \in \Hom_\cc(C^{\otimes n},A)$ of an object $n \in \N$ is given by the
convolution product:
\begin{equation*}
f \circ g = f \star g=\mu(f \otimes g)\Delta_{C^{\otimes n}}\co C^{\otimes n} \to A,
\end{equation*}
where $\Delta_{C^{\otimes n}}$ denotes the coproduct of the coalgebra $C^{\otimes n}$.

Note that the category $\conv_\cc(C,A)$ is a monoidal category: the monoidal product of two objects $m,n \in \N$ is given by
$m \otimes n=m+n$, the unit object is $0\in \N$, and the monoidal product of two morphisms $f\co m \to m$ and $g\co n \to n$
(where $m,n \in \N$) is the morphism $f \otimes g=\mu(f \otimes_\cc g)\co m+n \to m+n$.

\subsection{The category $\diag$} Let $\dd_0$ be the braided category freely generated by one object $*$ and the
following morphisms:
\begin{align*}
& \Delta\co  * \to * \otimes *, && \omega_+\co  * \otimes * \to \un, && \theta_+\co * \to \un, \\
& \varepsilon\co  * \to \un, && \omega_-\co * \otimes * \to \un, && \theta_-\co * \to \un,
\end{align*}
where $\un$ denotes the unit object of the monoidal product. Let $\dd$ be the quotient of the category $\dd_0$ by the
following relations:
\begin{align}
& (\id_* \otimes \Delta)\Delta=(\Delta \otimes \id_*)\Delta,\label{comult} \\
& (\id_* \otimes \varepsilon)\Delta=\id_*=(\varepsilon \otimes \id_*)\Delta. \label{counit}
\end{align}
The category $\dd$ is still braided (with induced braiding) and $(*,\Delta,\varepsilon)$ is a coalgebra in $\dd$. We define
the category $\diag$ to be the convolution category $\conv_{\dd}(*,\un)$, see Section~\ref{convcat}, where $\un$ is endowed
with the trivial algebra structure.

\subsection{The category $\diag^S$}\label{sectds} Let $\dd^S_0$ be the braided category freely generated by one object $*$ and the
following morphisms:
\begin{align*}
& \Delta\co  * \to * \otimes *, && \omega_+\co  * \otimes * \to \un, && \theta_+\co * \to \un, \\
& \varepsilon\co  * \to \un, && \omega_-\co * \otimes * \to \un, && \theta_-\co * \to \un,\\
& S\co * \to *, && S^{-1}\co  * \to *.
\end{align*}
Let $\dd^S$ be the quotient of the category $\dd^S_0$ by the relations \eqref{comult} and \eqref{counit}, and the following
relations:
\begin{align}
& SS^{-1}=\id_*=S^{-1}S,\label{antip1} \\
& \Delta S= (S \otimes S)\tau_{*,*} \Delta,\label{antip2} \\
& \varepsilon S=\varepsilon,\label{antip3} \\
& \theta_\pm S= \theta_\pm,\label{antip4} \\
& \omega_+ (S \otimes \id_*)=\omega_-=\omega_+ (\id_* \otimes S),\label{antip5} \\
& \omega_+ (S^{-1} \otimes \id_*)=\omega_- \tau_{*,*}=\omega_+ (\id_* \otimes S^{-1}),\label{antip6}
\end{align}
where $\tau_{*,*}\co * \otimes * \to * \otimes *$ denotes the braiding of the object $*$ with itself in~$\dd^S_0$.

The category $\dd^S$ is still braided (with induced braiding) and $(*,\Delta,\varepsilon)$ is a coalgebra in $\dd^S$. We
define the category $\diag^S$ to be the convolution category $\conv_{\dd^S}(*,\un)$, see Section~\ref{convcat}, where $\un$
is endowed with the trivial algebra structure.

\subsection{Relations between $\diag$ and $\diag^S$} The inclusion functor $\dd_0 \hookrightarrow \dd^S_0$ induces a functor
$\dd \to \dd^S$ and so a functor $\iota\co \diag \to \diag^S$. Note that $\iota$ is the identity on the objects.

\begin{thm}\label{diaganddiagS}
$\iota\co \diag \to \diag^S$ is an isomorphism of categories.
\end{thm}

\begin{proof}
Fix $n \in \N$. Set $C=\Hom_{\dd_0}(*^n,*)$ and $C^S = \Hom_{\dd^S_0}(*^n,*)$. We identify $C$ with its image under the
functor $\dd_0 \hookrightarrow \dd^S_0$, so that we have $C \subset C^S$. Let $\sim$ be the equivalence relation on $C$
defined by \eqref{comult}-\eqref{counit}, and let $\sim^S$ be the equivalence relation on $C^S$ defined by
\eqref{comult}-\eqref{antip6}. All we have to show is that the map $F:C/_\sim \to C^S/_{\sim^S}$ induced by the inclusion $C
\hookrightarrow C^S$ is bijective.

We will construct a map $C^S \to C$, $f \mapsto f'$, satisfying:
\begin{enumerate}
\renewcommand{\labelenumi}{{\rm (\alph{enumi})}}
\item $f' \sim^S f$ for all $f \in C^S$;
\item $f' = f$ for all $f \in C$;
\item $f \sim^S g \Longrightarrow f' \sim g'$ for all $f,g \in C^S$.
\end{enumerate}
If this holds, then $F$ is onto by (a), and $F$ is into by (b) and (c).

Now consider the following rewriting rules:
\begin{align*}
 & \Delta S^{\pm 1} \longrightarrow (S^{\pm 1} \otimes S^{\pm 1})\tau^{\pm 1}_{*,*} \Delta,
 && \epsilon S^{\pm 1} \longrightarrow \epsilon,\\
 & SS^{-1}\longrightarrow \id_*,
 && S^{-1}S \longrightarrow \id_*,\\
 & \theta_+ S^{\pm 1}\longrightarrow \theta_+,
 && \theta_- S^{\pm 1}\longrightarrow \theta_-,\\
 & \omega_+(S\otimes \id_*)\longrightarrow \omega_-,
 && \omega_+(\id_* \otimes S) \longrightarrow \omega_-,\\
 & \omega_+(S^{-1} \otimes \id_*)\longrightarrow \omega_- \tau_{*,*},
 &&  \omega_+ (\id_* \otimes S^{-1}) \longrightarrow \omega_- \tau_{*,*},\\
 & \omega_-(S \otimes \id_*)\longrightarrow \omega_+ \tau^{-1}_{*,*},
 &&  \omega_-(\id_* \otimes S)\longrightarrow  \omega_+ \tau^{-1}_{*,*}, \\
 & \omega_-(S^{-1} \otimes \id_*)\longrightarrow \omega_+,
 &&  \omega_-(\id_* \otimes S^{-1}) \longrightarrow \omega_+.
\end{align*}
For $f, g \in C^S$, we write $f \le g$ (resp.\@ $f \prec g$) if we can go from $f$ to $g$ by applying a finite number
(resp.\@ exactly one) of these rules. The rewriting system is noetherian. Indeed, the first rule decreases the number of
letters $S^{\pm 1}$ on the right of $\Delta$, and all other rules decrease the number of letters $S^{\pm 1}$. Therefore
$(C^S, \le)$ is a noetherian partially ordered set. Moreover the system is confluent: if $x \prec y$ and $x \prec z$, then
there exists $t$ such that $y, z \le t$. This can easily be checked by considering all cases when two rewriting rules can be
applied to the same element $f \in C^S$. By \cite{Ne}, for each $f \in C^S$ there exists a unique maximal element $f' \in
C^S$ such that $f \le f'$. The maximality condition means that no rewriting rule can be applied to $f'$. Such can only be the
case if no letter $S^{\pm 1}$ occurs in $f'$, that is, if $f' \in C$.

So we have constructed a map $f \mapsto f'$. Condition (b) is obvious. Condition~(a) results from the observation that
$\sim^S$ is the equivalence relation on $C^S$ generated by $\prec$ and \eqref{comult}-\eqref{counit}. Let us check Condition
(c). By construction, if $f \prec g$ then $f' = g'$. So, by the previous observation, it is enough to verify that if $g$ is
obtained from $f$ by a modification of the form (1) or (2), then $f' \sim g'$. This is easily shown by noetherian induction.
The only point here is to verify that, in the presence of a configuration of the type $(\Delta \otimes \id_*)\Delta S^{\pm
1}$, applying~\eqref{comult} followed by twice the first rewriting rule gives the same result as applying twice the first
rewriting rule followed by \eqref{comult}. Hence the theorem.
\end{proof}

\subsection{Hopf diagrams} 
\begin{figure}[ht!]\anchor{gensdiagS}\anchor{deltn}
   \begin{center}
       \subfigure[$\Delta \co * \to * \otimes *$]{\, \scaledraw{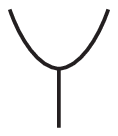} \,} \quad\;
       \subfigure[$\varepsilon \co * \to \un$]{\scaledraw{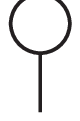}} \quad\;
       \subfigure[$S\co * \to *$]{\,\scaledraw{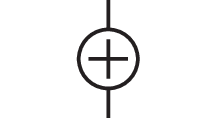}\,} \quad\,
       \subfigure[$S^{-1}\co * \to *$]{\;\scaledraw{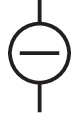}\;} \\
       \subfigure[$\Delta^{(n)} \co * \to *^{\otimes(n+1)}$]{\, \scaledraw{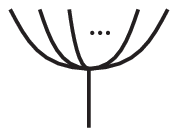} \, \label{deltn}}\quad\;
       \subfigure[$\tau_{*,*}\co * \otimes * \to * \otimes *$]{\,\scaledraw{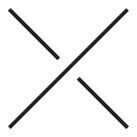}\,} \quad\;
       \subfigure[$\tau_{*,*}^{-1}\co * \otimes * \to * \otimes *$]{\scaledraw{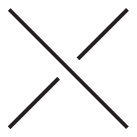}} \\
       \subfigure[$\omega_+\co * \otimes * \to \un$]{\,\scaledraw{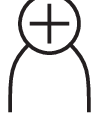}\;} \quad\;
       \subfigure[$\omega_-\co * \otimes * \to \un$]{\,\scaledraw{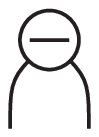}\,} \quad\;
       \subfigure[$\theta_+\co * \to \un$]{\,\scaledraw{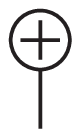}\,} \quad\;
       \subfigure[$\theta_-\co * \to \un$]{\;\scaledraw{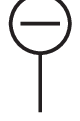}\,}
   \end{center}
     \caption{Generators of $\diag^S$}
     \label{gensdiagS}
\end{figure}

\begin{figure}[ht!]\anchor{relsdiagS}
   \begin{center}
       \scaleraisedraw{.403}{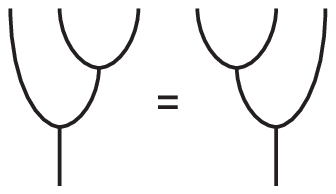}\, , \qquad
       \scaleraisedraw{.403}{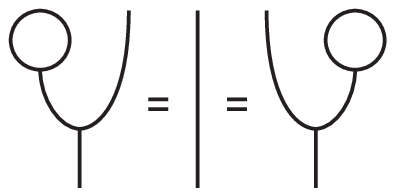}, \\[1.1em]
       \scaleraisedraw{.403}{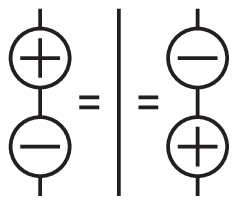}\; , \qquad
       \scaleraisedraw{.403}{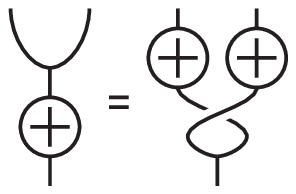} , \qquad
       \scaleraisedraw{.403}{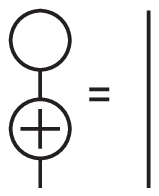}\; , \qquad
       \scaleraisedraw{.403}{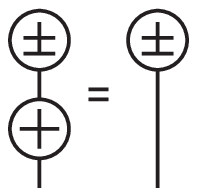} , \\[1.1em]
       \scaleraisedraw{.403}{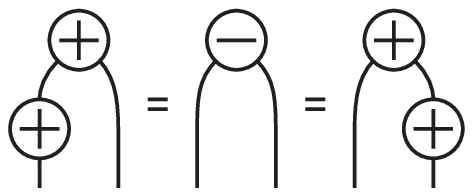}\; , \qquad
       \scaleraisedraw{.403}{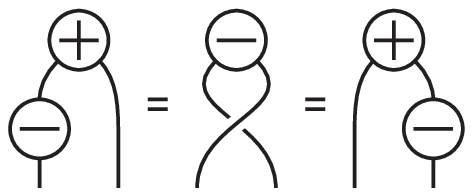}\; .
   \end{center}
     \caption{Relations in $\diag^S$}
     \label{relsdiagS}
\end{figure}

\begin{figure}[ht!]\small\anchor{convoprod}
   \begin{center}
       \psfrag{F}[cc][cc]{$D\star D'$}
       \psfrag{D}[cc][cc]{$D$}
       \psfrag{E}[cc][cc]{$D'$}
       \psfrag{o}[cc][cc]{$\circ$}
       \scaledraw{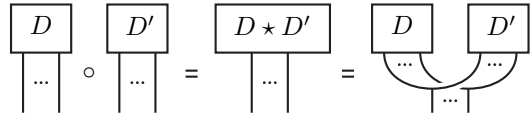}
   \end{center}
     \caption{Composition in $\diag^S$}
     \label{convoprod}
\end{figure}

By a \emph{Hopf diagram}, we shall mean a morphism of $\diag$ or $\diag^S$. Hopf diagrams can be
represented by plane diagrams: we draw one of their preimage in the braided category $\dd$ or $\dd^S$ by using Penrose
graphical calculus with the ascending convention (diagrams are read from bottom to top) as for instance in~\cite{Tur2}. We
depict the generators as in \figref{gensdiagS} except \figref{deltn} which depicts $\Delta^{(n)} \co
* \to *^{\otimes(n+1)}$ defined inductively by:
\begin{equation*}
\Delta^{(0)}=\id_*, \quad \Delta^{(1)}=\Delta \quad \text{and} \quad \Delta^{(n+1)}=(\Delta^{(n)}\otimes \id_*)\Delta.
\end{equation*}
The relations defining $\diag^S$ (except those concerning the braiding) are depicted in \figref{relsdiagS}. Recall that
the composition $D_1 \circ D_2=D_1\star D_2$ of two Hopf diagrams $D_1$ and $D_2$ is given by the convolution product, see
\figref{convoprod}.

\section{Ribbon string links and ribbon handles}\label{sect-hand}

\begin{figure}[ht!]\anchor{reidemeister}
   \begin{center}
      \subfigure[Type 0]{\scaledraw{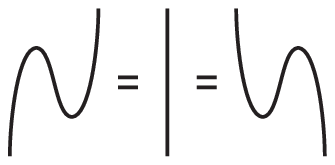}} \qquad\;
      \subfigure[Type 1$'$]{\scaledraw{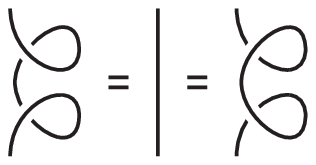}} \qquad\;
      \subfigure[Type 2]{\scaledraw{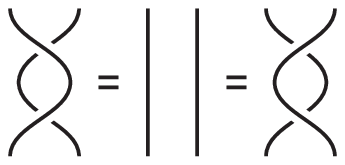}} \\
      \subfigure[Type 2$'$]{\scaledraw{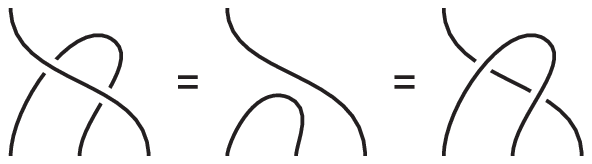}} \qquad\;
      \subfigure[Type 3]{\scaledraw{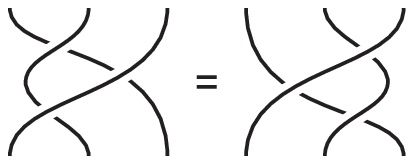}}
   \end{center}
     \caption{Ribbon Reidemeister moves}
     \label{reidemeister}
\end{figure}

As usual, we represent ribbon tangles (also called framed tangles) by thin plane diagrams (using blackboard framing). Recall
that two such diagrams represent the same isotopy class of a ribbon tangle if and only if one can be obtained from the other
by deformation and a finite sequence of ribbon Reidemeister moves depicted in \figref{reidemeister}.

\subsection{Ribbon string links}

Let $n$ be a non-negative integer. By a \emph{ribbon $n$-string link} we shall mean a ribbon $(n,n)$-tangle $T \subset
\mathbb{R}^2 \times [0,1]$  consisting of $n$ arc components, without any closed component, such that the $k$th arc ($1\leq k
\leq n$) joins the $k$th bottom endpoint to the $k$th top endpoint. Note that a ribbon string link is canonically oriented by
orienting each component from bottom to top.

We denote by $\rstl$ the category of ribbon string links. The objects of $\rstl$ are the non-negative integers. For two
non-negative integers $m$ and $n$, the set of morphisms from $m$ to $n$ is
\begin{equation*}
\Hom_{\rstl}(m,n)=\begin{cases} \emptyset & \text{if $m\neq n$,}\\ \rstl_n & \text{if $m=n$,} \end{cases}
\end{equation*}
where $\rstl_n$ denotes the set of (isotopy classes) of ribbon $n$-string links. The composition $T' \circ T$ of two ribbon
$n$-string links is given by stacking $T'$ on the top of $T$ (i.e., with ascending convention). Identities are the trivial
string links.

Note that the category $\rstl$ is a monoidal category: $m \otimes n=m+n$ on objects and the monoidal product $T \otimes T'$
of two ribbon string links $T$ and $T'$ is the ribbon string link obtained by juxtaposing $T$ on the left of $T'$ (see,
e.g., \cite{Tur2}).

\subsection{Ribbon handles}\label{sectrhand} Let $n$ be a non-negative integer. By a \emph{ribbon $n$-handle} we shall mean a
ribbon $(2n,0)$-tangle $T \subset \mathbb{R}^2 \times [0,1]$  consisting of $n$ arc components, without any closed component,
such that the $k$\trait th arc ($1\leq k \leq n$) joins the $(2k-1)$\trait th bottom endpoint to the $2k$\trait th bottom
endpoint. Note that a ribbon handle is canonically oriented by orienting each component upwards near its right bottom input.

\begin{figure}[ht!]
\anchor{comporhand}\anchor{identrhand}\anchor{tensorhand}
   \begin{center}
       \subfigure[$T$]{\psfrag{D}[cc][cc]{\small$T$}
                       \scaledraw{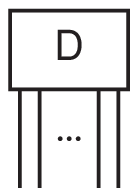}}
       \qquad
       \subfigure[$T'$]{\psfrag{E}[cc][cc]{$T'$}
                        \scaledraw{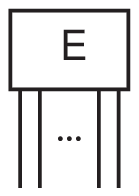}}
       \qquad
       \subfigure[$T \circ T'$]{\psfrag{D}[cc][cc]{\small$T$}
                                \psfrag{E}[cc][cc]{\small$T'$}
                                \psfrag{o}[cc][cc]{\small$\circ$}
                                \scaledraw{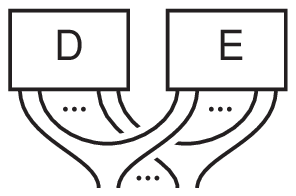}
                                \label{comporhand}}
       \qquad
       \subfigure[$\id$]{\scaledraw{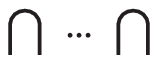}
                         \label{identrhand}}
       \qquad
       \subfigure[$T \otimes T'$]{\psfrag{D}[cc][cc]{\small$T$}
                                \psfrag{E}[cc][cc]{\small$T'$}
                                \scaledraw{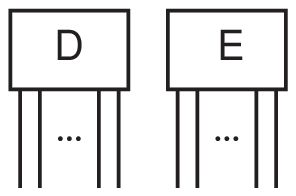}
                                \label{tensorhand}}
   \end{center}
     \caption{Composition, identity, and monoidal product in  $\rhand$}
\end{figure}

We denote by $\rhand$ the category of ribbon handles. The objects of $\rhand$ are the non-negative integers. For two
non-negative integers $m$ and $n$, the set of morphisms from $m$ to $n$ is
\begin{equation*}
\Hom_{\rhand}(m,n)=\begin{cases} \emptyset & \text{if $m\neq n$,}\\ \rhand_n & \text{if $m=n$,} \end{cases}
\end{equation*}
where $\rhand_n$ denotes the set of (isotopy classes) of ribbon $n$-handles. The composition of two ribbon $n$-handles $T$
and $T'$ is the ribbon $n$-handle defined in \figref{comporhand}. The identity for this composition consists in $n$
caps, see \figref{identrhand}.

Note that the category $\rhand$ is a monoidal category: $m \otimes n=m+n$ on objects and the monoidal product $T \otimes T'$
of two ribbon handles $T$ and $T'$ is the ribbon handle obtained by juxtaposing $T$ on the left of $T'$, see
\figref{tensorhand}.

\subsection{An isomorphism between ribbon handles and string links}\label{FandG}

\begin{figure}[ht!]\anchor{rstl2rhand}\anchor{rhand2rstl}
   \begin{center}
       \subfigure[$S$]{\psfrag{E}[cc][cc]{\small$S$}
                          \scaledraw{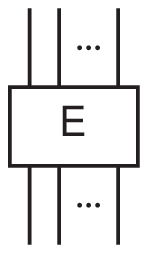}}
       \qquad \qquad
       \subfigure[$F(S)$]{\psfrag{E}[cc][cc]{\small$S$}
                          \scaledraw{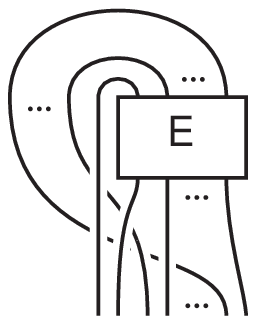}
                          \label{rstl2rhand}}
       \qquad \qquad
       \subfigure[$T$]{\psfrag{T}[cc][cc]{\small$T$}
                          \scaledraw{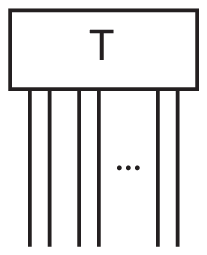}}
       \qquad \qquad
       \subfigure[$G(T)$]{\psfrag{T}[cc][cc]{\small$T$}
                          \scaledraw{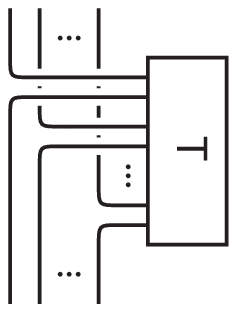}
                          \label{rhand2rstl}}
   \end{center}
     \caption{Definition of the functors $F$ and $G$}
\end{figure}

Let us construct functors $F\co \rstl \to \rhand$ and $G\co \rhand \to \rstl$ as follows. On objects, set $F(n)=n$ and
$G(n)=n$ for any non-negative integer $n$. For any ribbon $n$-string link $S$, let $F(S)$ be the ribbon $n$-handle defined
in \figref{rstl2rhand}. For any ribbon $n$-handle $T$, let $G(T)$ be the ribbon $n$-string link defined in
\figref{rhand2rstl}.

\begin{prop}\label{propFandG}
The functors $F\co \rstl \to \rhand$ and $G\co \rhand \to \rstl$ are mutually inverse monoidal functors.
\end{prop}
\begin{proof}Straightforward.
\end{proof}

\section{From Hopf diagrams to ribbon string links}\label{sect-hopf2hand}

\subsection{From Hopf diagrams to ribbon handles}
\begin{figure}[ht!]\anchor{Hopf2Hand}
   \begin{center}
      \psfrag{,}[Bl][Bl]{,}
      \psfrag{.}[Bl][Bl]{.}
      \scaledraw{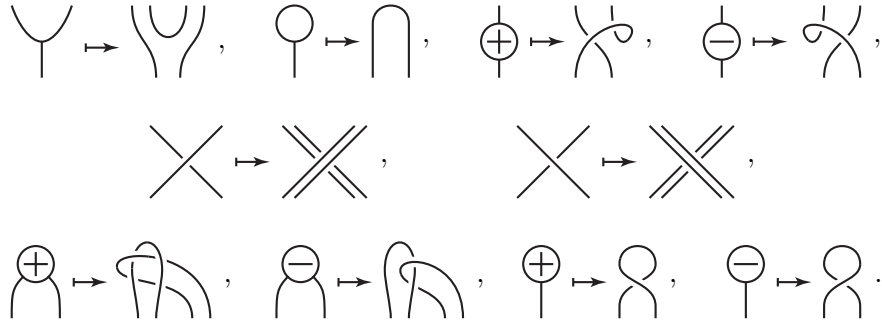}
   \end{center}
     \caption{Rules for defining $\phi$}
     \label{Hopf2Hand}
\end{figure}

\begin{figure}[ht!]\anchor{ExHopf2Hand}
   \begin{center}
      $D=$  \scaleraisedraw{.37}{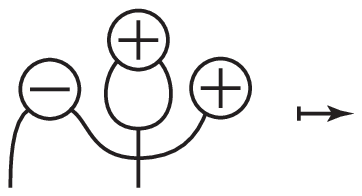}  \, $\phi_D= \,$
      \psfrag{-}[cc][cc]{\scalebox{1.25}{$\sim$}}\scaleraisedraw{.37}{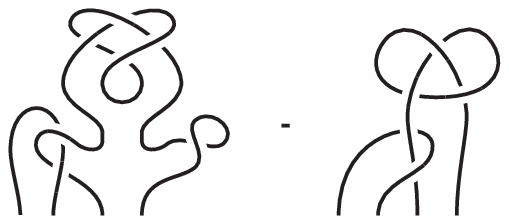}
   \end{center}
     \nocolon\caption{}
     \label{ExHopf2Hand}
\end{figure}

Let us define a functor $\phi$ from the category $\diag^S$ of Hopf diagrams
to the category $\rhand$ of ribbon handles. For any non-negative integer $n$, we set $\phi(n)=n$. Given a Hopf diagram $D$,
we construct a diagram of ribbon handle $\phi_D$ by using the rules of \figref{Hopf2Hand} and the stacking product (with
ascending convention). See \figref{ExHopf2Hand} for an example. Then let $\phi(D)$ be the isotopy class of the ribbon
handles defined by $\phi_D$.

\begin{prop}\label{Hopf2handthm}
The functor $\phi\co \diag^S \to \rhand$ is well-defined and mo\-noidal.
\end{prop}

We will prove in Section~\ref{sect-hand2hopf} (see Corollary~\ref{rhand2diag}) that $\phi$ is surjective.

\begin{figure}[ht!]\anchor{Phiwelldef}
   \begin{center}
      \psfrag{,}[Bl][Bl]{,}
      \psfrag{.}[Bl][Bl]{.}
      \psfrag{-}[B][B]{$\sim$}
      \scaledraw{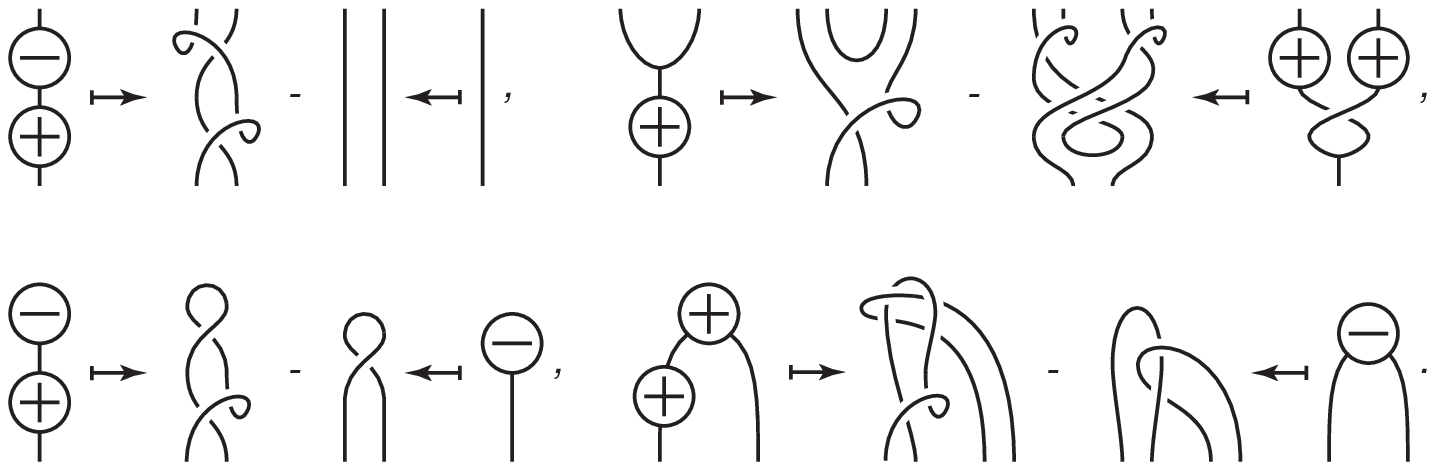}
   \end{center}
     \nocolon\caption{}
     \label{Phiwelldef}
\end{figure}

\begin{proof}
Let us first verify that $\phi$ is well-defined. We only have to verify that both sides of the relations defining $\diag^S$
are transformed by the rules of \figref{Hopf2Hand} to isotopic tangles. Examples of such verifications are depicted in
\figref{Phiwelldef}. The fact that $\phi$ is a monoidal functor comes from the definitions of composition and monoidal
product in $\diag^S$ and $\rhand$.
\end{proof}

\subsection{From Hopf diagrams to ribbon string links} 
\begin{figure}[ht!]\anchor{functpsi}\anchor{functpsiEx}
   \begin{center}
       \subfigure[]{\psfrag{T}[cc][cc]{\small$\phi_D$}
       \psfrag{p}[cc][cc]{\small$\psi_D$}
       \scaledraw{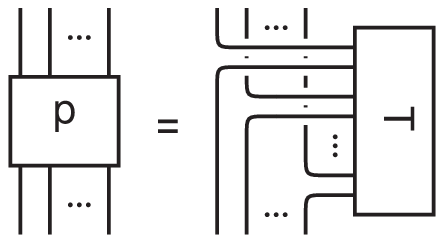}\label{functpsi}}\qquad \qquad
       \subfigure[]{\psfrag{e}[Bl][Bl]{$\rightsquigarrow$}
                            \psfrag{-}[c][c]{$\sim$}
                            \scaledraw{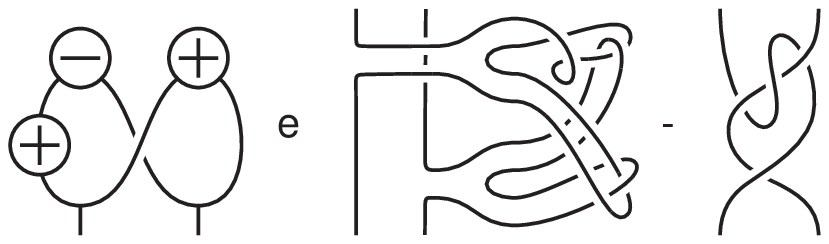} \label{functpsiEx}}
   \end{center}
   \nocolon\caption{}
\end{figure}

Let us define a functor $\psi \co  \diag^S \to \rstl$ as follows. For any
non-negative integer $n$, let $\psi(n)=n$. If $D$ is a Hopf diagram, we construct a diagram of ribbon string link $\psi_D$
as in \figref{functpsi}, where $\phi_D$ is defined as above. For an example, see \figref{functpsiEx}. Then we let
$\psi(D)$ be the isotopy class of the ribbon string link defined by $\psi_D$.

\begin{cor}
The functor $\psi\co \diag^S \to \rstl$ is well-defined and surjective. Moreover $F \circ \psi=\phi$ and $G \circ
\phi=\psi$, where $F$ and $G$ are the functors of Section~\ref{FandG}.
\end{cor}
\begin{proof}
This is an immediate consequence of Propositions~\ref{propFandG} and \ref{Hopf2handthm}.
\end{proof}

\section{From ribbon string links to Hopf diagrams}\label{sect-hand2hopf}
In this section, we construct a refined version $\qdiagS$ of the category of Hopf diagrams and we prove the main theorem.

\subsection{More relations on Hopf diagrams}\label{reladdhopfdiag}
\begin{figure}[ht!]\anchor{relssup}
   \begin{center}
       \psfrag{,}[Bl][Bl]{,}
       \psfrag{.}[Bl][Bl]{.}
       \scaledraw{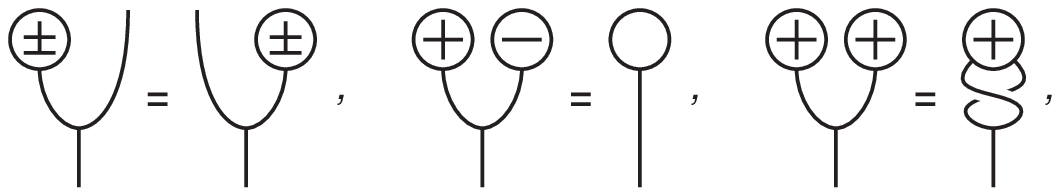}\\[1em]
       \scaledraw{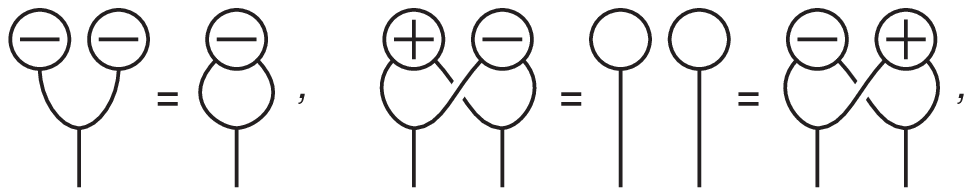} \\[1em]
       \scaledraw{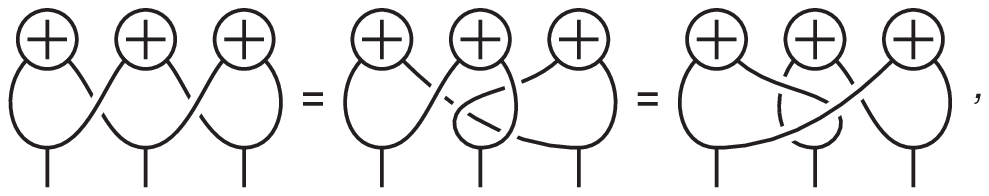} \\[1em]
       \scaledraw{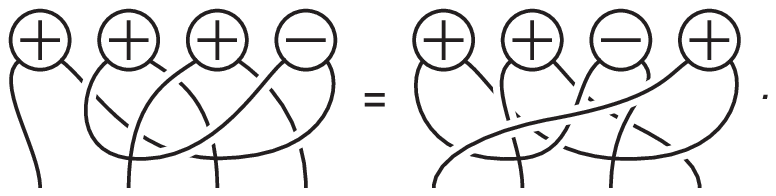}
   \end{center}
     \nocolon\caption{}
     \label{relssup}
\end{figure}

Recall the braided category $\dd^S$ of Section~\ref{sectds}. Let $\qddS$ be the quotient of $\dd^S$ by the following
relations:
\begin{align}
& (\id_* \otimes \theta_\pm)\Delta=(\theta_\pm \otimes \id_*)\Delta,\label{relsuptwist1} \\
& (\theta_+ \otimes \theta_-)\Delta=\varepsilon,\label{relsuptwist2} \\
& (\theta_+ \otimes \theta_+)\Delta=\omega_+ \tau_{*,*}^{-1}\Delta,\label{relsuptwist3} \\
& (\theta_- \otimes \theta_-)\Delta=\omega_- \Delta,\label{relsuptwist4} \\
& \omega_+ \star \omega_-=\varepsilon \otimes \varepsilon=\omega_- \star \omega_+, \label{relsupRII}\\
& \omega_{+12} \star \omega_{+13} \star \omega_{+23}= \omega_{+13} \star \omega_{+23} \star \omega_{+12}
= \omega_{+23} \star \omega_{+12} \star \omega_{+13},\label{relsup3T}\\
& \omega_{+13} \star \omega_{+23} \star \omega_{+24} \star \omega_{-23}
  =\omega_{+23} \star \omega_{+24} \star \omega_{-23} \star \omega_{+13}, \label{relsup4T}
\end{align}
where $\star$ denotes the convolution product of Hopf diagrams. Graphically, these relations can be depicted as in
\figref{relssup}. The category $\qddS$ is still braided (with induced braiding) and $(*,\Delta,\varepsilon)$ is a
coalgebra in $\qddS$. Then let $\qdiagS$ be the convolution category $\conv_{\qddS}(*,\un)$, see Section~\ref{convcat}, where
$\un$ is endowed with the trivial algebra structure. By abuse,  we will still call \emph{Hopf diagrams} the morphisms of
$\qdiagS$.

Note that the category $\qdiagS$ can also be viewed as the quotient of the category $\diag^S$ by
Relations~\eqref{relsuptwist1}-\eqref{relsup4T}. Let $\pi\co \diag^S \to \qdiagS$ be the projection functor.

\begin{prop}
The functors $\phi$ and $\psi$ factorize through $\pi$ to functors $\overline{\phi}$ and $\overline{\psi}$, respectively,
so that $\overline{\phi}=F \circ \overline{\psi}$.
\end{prop}
\begin{proof}
We have to verify both sides of Relations~\eqref{relsuptwist1}-\eqref{relsup4T} are transformed by the rules of
\figref{Hopf2Hand} to isotopic tangles. Examples of such verifications are depicted in \figref{Moreker}. Note that
in the first picture, we used the fact that the $(2i-1)$-th and $2i$-th inputs of a ribbon handle are connected.
\begin{figure}[ht!]\anchor{demquotphi1bis}
   \begin{center}
      \psfrag{,}[Bl][Bl]{,}
      \psfrag{.}[Bl][Bl]{.}
      \psfrag{-}[B][B]{$\sim$}
      \scaledraw{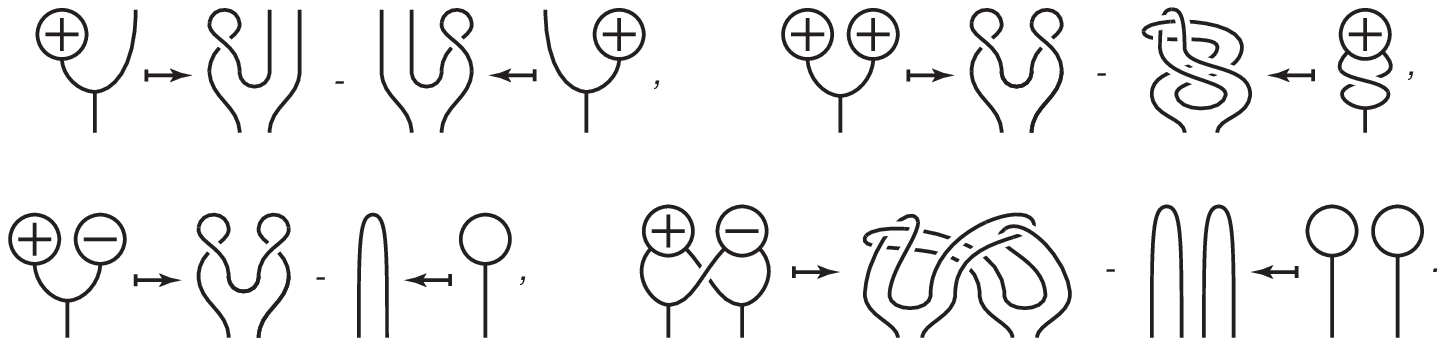}
   \end{center}
     \nocolon\caption{}
     \label{Moreker}
\end{figure}
\end{proof}

\subsection{From ribbon pure braids to Hopf diagrams}\label{psirpb}
Let $\RPB$ the subcategory of $\rstl$ made of ribbon pure braids. The objects of $\RPB$ are the non-negative integers. For
two non-negative integers $m$ and $n$, the set of morphisms from $m$ to $n$ is
\begin{equation*}
\Hom_{\RPB}(m,n)=\begin{cases} \emptyset & \text{if $m\neq n$,}\\ \RPB_n & \text{if $m=n$,} \end{cases}
\end{equation*}
where $\RPB_n \subset \rstl_n$ denotes the set of (isotopy classes) of ribbon pure $n$-braids. Note that $\RPB$ is a
monoidal subcategory of $\rstl$.

Recall we have a canonical group isomorphism:
\begin{equation*}
(u,t_1, \dots, t_n) \co  \RPB_n \stackrel{\sim}{\longrightarrow} \PB_n \times \mathbb{Z}^n,
\end{equation*}
where $\PB_n$ denotes the group of pure $n$-braids, $u\co \RPB_n \to \PB_n$ is the forgetful morphism, and $t_i$ the
self-linking number of the $i$-th component. Hence, using a presentation of $\PB_n$ by generators and relations due to
Markov \cite{Mark1}, we get that $\RPB_n$ is generated by $t_k$ ($1\leq k \leq n$) and $\sigma_{i,j}$ ($1 \leq i < j \leq
n$) subject to the following relations:
\begin{align}
& t_k t_l = t_l t_k \quad \text{for any $k$, $l$;} \label{Markoff1}\\
& t_k \sigma_{i,j} = \sigma_{i,j} t_k \quad \text{for any $i<j$ and $k$;} \label{Markoff2}\\
& \sigma_{i,j}\sigma_{k,l}=\sigma_{k,l}\sigma_{i,j} \quad \text{for any $i < j< k< l$ or any $i<k<j<l$;}\label{Markoff3}\\
& \sigma_{i,j}\sigma_{i,k}\sigma_{j,k}=\sigma_{i,k}\sigma_{j,k}\sigma_{i,j}=\sigma_{j,k}\sigma_{i,j}\sigma_{i,k} \quad
\text{for any $i<j<k$;}\label{Markoff4}\\
& \sigma_{i,k}\sigma_{j,k}\sigma_{j,l}\sigma^{-1}_{j,k}=\sigma_{j,k}\sigma_{j,l}\sigma^{-1}_{j,k}\sigma_{i,k} \quad
\text{for any $i<j<k<l$.}\label{Markoff5}
\end{align}
Graphically, the generators may be represented as:
\begin{equation*}
      \sigma_{i,j}=
      \psfrag{1}[Bl][Bl]{{\scriptsize $1$}} \psfrag{i}[Bl][Bl]{{\scriptsize $i$}} \psfrag{j}[Bl][Bl]{{\scriptsize $j$}}
      \psfrag{n}[Bl][Bl]{{\scriptsize $n$}} \scaleraisedraw{.578}{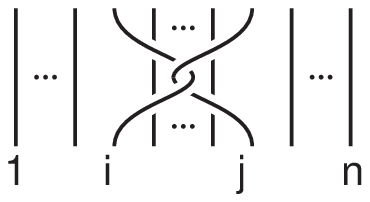}
      \qquad \text{and} \qquad t_{k}=
      \psfrag{1}[Bl][Bl]{{\scriptsize $1$}} \psfrag{k}[Bl][Bl]{{\scriptsize $k$}}
      \psfrag{n}[Bl][Bl]{{\scriptsize $n$}} \scaleraisedraw{.578}{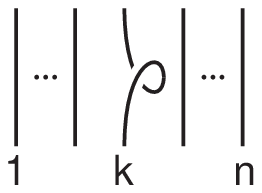} .
\end{equation*}

Let us define a functor $\Psi_0\co  \RPB \to \qdiagS$ as follows: for any non-negative integer $n$ set $\Psi(n)=n$. For $1
\leq i < j\leq n$ and $1 \leq k \leq n$, set
\begin{equation*}
\Psi_0(\sigma_{i,j}^{\pm 1})=\Sigma_{i,j}^{\pm 1} \quad \text{and} \quad \Psi_0(t_{k}^{\pm 1})=\Omega_k^{\pm 1},
\end{equation*}
where
\begin{equation}\label{gammaomega}
      \Sigma_{i,j}^{\pm 1}=
      \psfrag{1}[Bl][Bl]{{\scriptsize $1$}} \psfrag{i}[Bl][Bl]{{\scriptsize $i$}} \psfrag{j}[Bl][Bl]{{\scriptsize $j$}}
      \psfrag{n}[Bl][Bl]{{\scriptsize $n$}} \scaleraisedraw{.4}{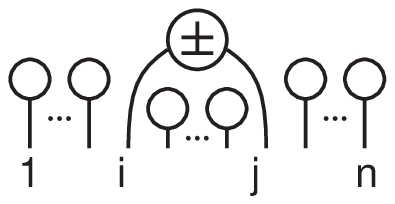}
      \qquad \text{and} \qquad \Omega_k^{\pm 1}=
      \psfrag{1}[Bl][Bl]{{\scriptsize $1$}} \psfrag{k}[Bl][Bl]{{\scriptsize $k$}}
      \psfrag{n}[Bl][Bl]{{\scriptsize $n$}} \scaleraisedraw{.4}{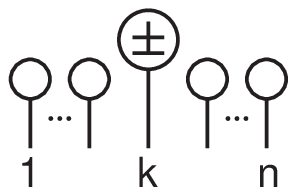} .
\end{equation}

\begin{lem}\label{proppure2diag}
The functor $\Psi_0\co \RPB \to \qdiagS$ is well defined, monoidal, and is such that $\overline{\psi} \circ \Psi_0(P)=P$ for
all ribbon pure braid $P$.
\end{lem}
\begin{proof}
Firstly, from Relation \eqref{relsupRII} (resp.\@ Relations \eqref{relsuptwist1} and \eqref{relsuptwist2}), we see that the
Hopf diagrams $\Sigma_{i,j}^{- 1}$ and $\Sigma_{i,j}$ (resp.\@ $\Omega_k^{- 1}$ and $\Omega_k$) are inverse each other.

Secondly Relations \eqref{Markoff1}-\eqref{Markoff5} hold in $\qdiagS$, where we replace $\sigma_{i,j}$ and $t_k$ with
$\Sigma_{i,j}$ and $\Omega_k$ respectively. Indeed Relations \eqref{Markoff1} and \eqref{Markoff2} follow from
\eqref{relsuptwist1}. Relation \eqref{Markoff3} follows from \eqref{counit}. Relations \eqref{Markoff4} and \eqref{Markoff5}
correspond to \eqref{relsup3T} and \eqref{relsup4T} respectively.

Finally the isotopies depicted in \figref{demoprop2diagfig} show that $\overline{\psi} \circ
\Psi_0(\sigma_{i,j})=\sigma_{i,j}$ and $\overline{\psi} \circ \Psi_0(t_k)=t_k$. Hence $\overline{\psi} \circ \Psi_0(P)=P$
for all ribbon pure braid $P$.
\begin{figure}[ht!]\anchor{demoprop2diagfig}
   \begin{center}
                     \psfrag{1}[Bl][Bl][.8]{$1$}
                     \psfrag{i}[Bl][Bl][.8]{$i$}
                     \psfrag{j}[Bl][Bl][.8]{$j$}
                     \psfrag{n}[Bl][Bl][.8]{$n$}
                     \psfrag{-}[B][B]{$\sim$}
                     \scaleraisedraw{.54}{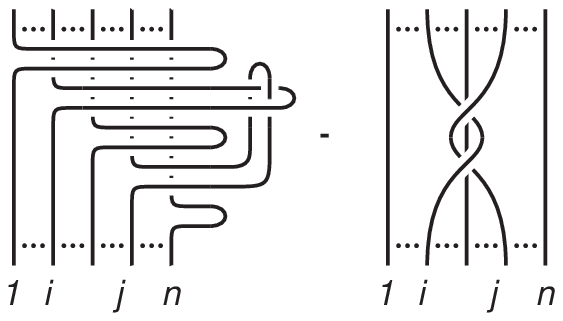}
                     \qquad \qquad
                     \psfrag{1}[Bl][Bl][.8]{$1$}
                     \psfrag{k}[Bl][Bl][.8]{$k$}
                     \psfrag{n}[Bl][Bl][.8]{$n$}
                     \psfrag{-}[B][B]{$\sim$}
                     \scaleraisedraw{.54}{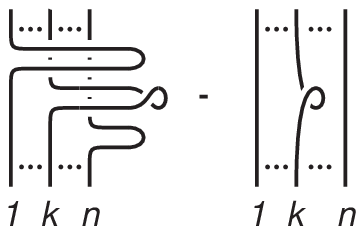}
   \end{center}
   \nocolon\caption{} \label{demoprop2diagfig}
\end{figure}
\end{proof}

\subsection{Contractions}\label{contracts}
\begin{figure}[ht!]\anchor{contractPB}\anchor{contractHD}
   \begin{center}
       \subfigure[]{\psfrag{P}[cl][cl]{\small$T$}
                            \psfrag{g}[cr][cr]{$c_i(T)=$}
                            \psfrag{a}[cc][cc]{$\overbrace{\hspace{7.1mm}}^{i-2}$}
                            \psfrag{c}[cc][cc]{$\overbrace{\hspace{7.1mm}}^{n-i-1}$}
                            \psfrag{o}[cc][cc]{$\underbrace{\hspace{7.1mm}}_{i-2}$}
                            \psfrag{e}[cc][cc]{$\underbrace{\hspace{7.1mm}}_{n-i-1}$}
                            \scaledraw{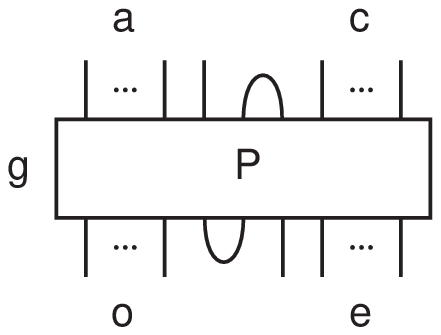}
                            \label{contractPB}}
       \qquad \qquad
       \subfigure[]{\psfrag{D}[cl][cl]{\small$D$}
                             \psfrag{g}[cr][cr]{$C_i(D)=$}
                            \psfrag{o}[cc][cc]{$\underbrace{\hspace{7.1mm}}_{i-2}$}
                            \psfrag{e}[cc][cc]{$\underbrace{\hspace{7.1mm}}_{n-i-1}$}
                            \scaledraw{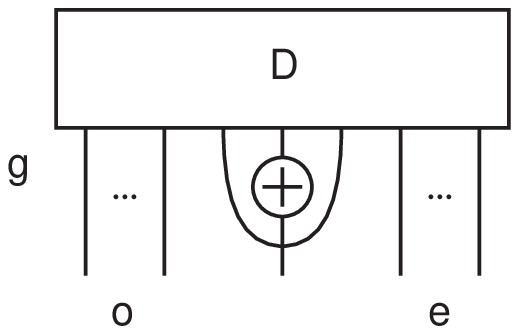}
                            \label{contractHD}}
   \end{center}
   \caption{Contractions of ribbon string links and Hopf diagrams}
\end{figure}

Let $n\geq 3$ and $1 < i < n$. For a ribbon $n$-string link $T$, we define the \emph{$i$-th contraction} of $T$ to be the
ribbon $(n-2)$-string link $c_i(T)$ defined as in \figref{contractPB}. For a Hopf diagram $D$ with $n$ inputs, we define
the \emph{$i$-th contraction} of $D$ to be the Hopf diagram with $(n-2)$ inputs $C_i(D)$ defined as in
\figref{contractHD}.

\begin{lem}\label{lemcontr}
Let $n\geq 3$ and $1 < i < n$. For any Hopf diagram $D$ with $n$ inputs, we have $c_i(\psi(D))=\psi(C_i(D))$.
\end{lem}
\begin{proof}
This follows from the following equalities:
\begin{equation*}
      \psi(C_i(D))=
      \psfrag{T}[cc][cc]{\small$\phi(D)$}\psfrag{-}[B][B]{$\sim$}
       \psfrag{E}[cc][cc]{\small$\phi(C_i(D))$} \scaleraisedraw{.47}{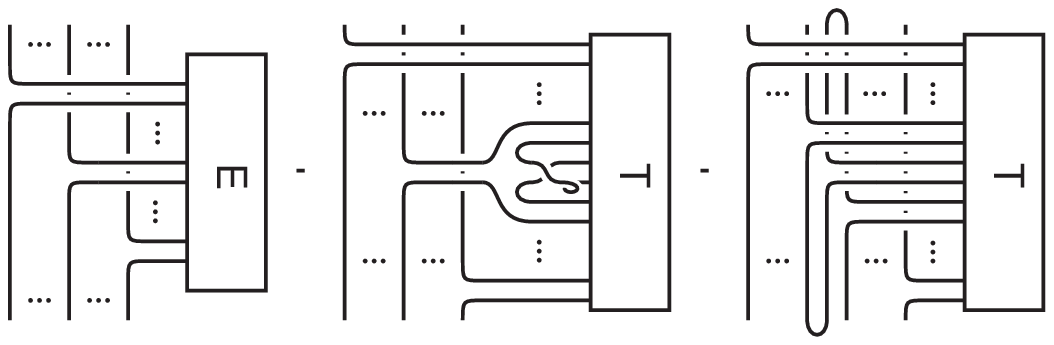} =c_i(\psi(D)),
\end{equation*}
where we used that the $(2i-1)$-th and $2i$-th inputs of $\phi(D)$ are connected (since it is a ribbon handle).
\end{proof}

\begin{lem}\label{contraccomm}
We have $c_ic_j= c_j c_{i+2}$ and $C_iC_j= C_j C_{i+2}$ for any $i \geq j$.
\end{lem}
\begin{proof}
This results directly from the definitions of the contraction operators.
\end{proof}

\subsection{From ribbon string links to Hopf diagram}\label{algorstl2diag}
In this section, we extend the monoidal functor $\Psi_0\co \RPB \to \qdiagS$ to a monoidal functor $\Psi\co \rstl \to
\qdiagS$. The construction of this extension follows, broadly speaking, the same pattern as the proof of \cite[Theorem
3]{Brug3}. The point is to see that a ribbon string link can be obtained from a ribbon pure braid by a sequence of
contractions. This will at least show that $\Psi$ extends uniquely and suggest a construction for it. We then must check the
coherence of this construction, that is, its independence from the choices we made.

\begin{figure}[ht!]\anchor{pullmax1}\anchor{pullmax2}\anchor{pullmax3}
   \begin{center}
       \subfigure[$\Gamma$]{\psfrag{E}[cl][cl]{\small$U$}
                    \psfrag{D}{\small$V$}
                    \psfrag{m}[cc][cc]{{\scriptsize $m$}}
                    \scaledraw{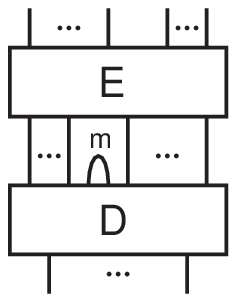}
                    \label{pullmax1}}
       \qquad \qquad
       \subfigure[$U \cup \ell$]{\psfrag{E}[cl][cl]{\small$U$}
                    \psfrag{l}{\small$\ell$}
                    \psfrag{1}[Bl][Bl]{{\scriptsize $1$}}
                    \psfrag{i}[Bl][Bl]{{\scriptsize $i$}}
                    \psfrag{j}[Bl][Bl]{{\scriptsize $j$}}
                    \psfrag{n}[Bl][Bl]{{\scriptsize $n$}}
                    \scaledraw{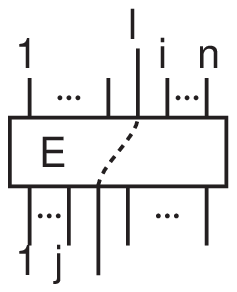}
                    \label{pullmax3}}
       \qquad \qquad
       \subfigure[$\Gamma_\ell$]{\psfrag{E}[cl][cl]{\small$U_\ell$}
                    \psfrag{D}{$V$}
                    \scaledraw{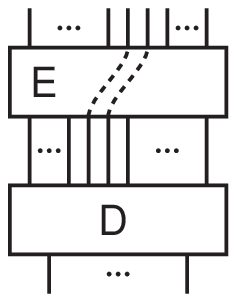}
                    \label{pullmax2}}
   \end{center}
   \caption{Pulling a max to the top line}
\end{figure}

The main trick we use consists in ``pulling a max to the top line''. Let $\Gamma$ be a tangle diagram with a local max $m$,
with $n$ outputs. We may write $\Gamma$ as in \figref{pullmax1}, where $U$, $V$ are tangle diagrams. Let $i$ be an
integer, $1 \leq i \leq n+1$. Let $j$ be the number of strands to the left of $m$ on the same horizontal line. Let $U \cup
\ell$ be a tangle diagram obtained from $U$ by inserting a new component $\ell$ going from a point between the $j$-th and
$(j+1)$-th inputs of $U$ to a point between the $(i-1)$-th and $i$-th outputs of $U$, see \figref{pullmax3}. We assume
also that $\ell$ has no local extremum. Let $U_\ell=\Delta_\ell(U \cup \ell)$ be the tangle diagram obtained from $U \cup
\ell$ by doubling $\ell$. Set $\Gamma_\ell=U_\ell V$, see \figref{pullmax2}. We say that \emph{$\Gamma_\ell$ is obtained
from $\Gamma$ by pulling $m$ to the top in the $i$-th position (along the path $\ell$)}. Likewise, one defines the action of
\emph{pulling a local min to the bottom.}

We define $\Psi\co \rstl \to \qdiagS$ as follows: on objects $n \in \N$, set $\Psi(n)=n$. Let~$n$ be a non-negative integer
and $T$ be a ribbon $n$-string link. Consider a diagram of~$T$. For each local extremum pointing to the right (once the
strands canonically oriented from bottom to top), modify the diagram using the following rule:
\begin{equation}\label{renderhanded}
\scaleraisedraw{.4}{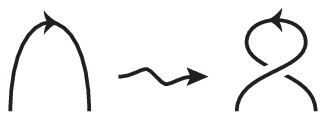} \quad \text{ or } \quad \scaleraisedraw{.4}{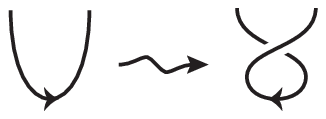}.
\end{equation}
This leads to a diagram $\Gamma$ which is  \emph{left handed}, that is, with all local extrema pointing to the left. Pulling
all local max to the top and all local min to the bottom, we obtain a diagram of a pure braid. Here is an algorithm. Denote
by $m_i$ the number of local max (which is equal to the number of local min) on the $i$\trait th component of $\Gamma$. Let
$m(\Gamma)=m_1 + \dots + m_n$ be the number of local max of $\Gamma$. If $m(\Gamma)=0$, we are already done. Otherwise,
chose $i$ maximal so that $m_i> 0$. Let $m$ be the first max and $m'$ be the first min you meet on the $i$-th component,
going from bottom to top. Pull $m$ to the top, in the $(i+1)$-th position, and $m'$ to the bottom, in the $i$-th position.
Let $\Gamma'$ be the diagram so constructed. Then $\Gamma'$ is a string link diagram, with $m(\Gamma')=m(\Gamma)-1$. Let us
denote by $\{\Gamma\}$ the ribbon string link defined by the diagram $\Gamma$. Then $\{\Gamma\}= c_{i+1} \{\Gamma'\}$, where
$c_{i+1}$ is the $(i+1)$-th contraction as in Section~\ref{contracts}. Repeating $m(\Gamma)$ times this transformation
yields a pure braid diagram $P$ with $n+ 2m(\Gamma)$ strands, and we have $\{\Gamma\}= c_{j_{m(\Gamma)}} \cdots c_{j_1}
\{P\}$ where  $1 \leq j_1 \leq \dots \leq j_{m(\Gamma)} \leq n$, and $j_k$ takes $m_i$ times the value $i+1$. Note that
\begin{equation}\label{formalgo}
T=(t_1^{\alpha_1} \cdots t_n^{\alpha_n}) \, c_{j_{m(\Gamma)}} \cdots c_{j_1} \{P\},
\end{equation}
where $\alpha_i$ is the number of modifications \eqref{renderhanded} made on the $i$-th component and $t_i\in \RPB_n$ is as
in Section~\ref{psirpb}. Finally, as suggested by Lemmas~\ref{proppure2diag} and \ref{lemcontr}, set:
\begin{equation*}
\Psi(T)= (\Omega_1^{\alpha_1}\cdots \Omega_n^{\alpha_n}) \, C_{j_{m(\Gamma)}} \cdots C_{j_1} \Psi_0 (\{P\}),
\end{equation*}
where $\Psi_0\co \RPB \to \qdiagS$ is the functor of Lemma~\ref{proppure2diag}, the $\Omega_j=\Psi_0(t_j)$ are the Hopf
diagrams of Section~\ref{psirpb}, and the $C_j$ are the contractions on Hopf diagrams defined in Section~\ref{contracts}.

A component of a ribbon string link is said to be \emph{trivial} if it is trivial as an unframed long knot (up to isotopy).
Note that each contraction $c_{j_k}$ in \eqref{formalgo} is performed on a ribbon string link with trivial $j_k$-th
component.

\begin{thm}\label{rstl2diag}
The functor $\Psi\co \rstl \to \qdiagS$ is well-defined, monoidal, and satisfies $\overline{\psi}\circ\Psi=1_{\rstl}$.
Moreover, $\Psi$ is the unique functor from $\rstl$ to $\qdiagS$ satisfying:
\begin{enumerate}
\renewcommand{\labelenumi}{{\rm (\alph{enumi})}}
  \item $\Psi(P)=\Psi_0(P)$ for any ribbon pure braid $P$;
  \item $\Psi(c_i(T))=C_i\Psi(T)$ for any $1<i < n$ and any ribbon $n$-string link $T$ with trivial $i$-th component.
\end{enumerate}
\end{thm}
Remark in particular that $\Psi$ is injective and $\overline{\psi}$ (and so $\psi$) is surjective.

We prove the theorem in Section~\ref{proofrstl2diag}.\\

Set $\Phi=\Psi \circ G \co  \rhand \to \qdiagS$, where $G\co \rhand \to \rstl$ is the monoidal isomorphism defined in
Section~\ref{FandG}. From Proposition~\ref{propFandG} and Theorem~\ref{rstl2diag}, we immediately deduce that:

\begin{cor}\label{rhand2diag}
The functor $\Phi\co \rhand \to \qdiagS$ is monoidal and satisfies $\overline{\phi}\circ\Phi=1_{\rhand}$.
\end{cor}
\noindent Note in particular that $\Phi$ is injective and that $\overline{\phi}$ (and so $\phi$) is surjective.

\subsection{Summary}
The previous results may be summarized in the commutativity of the following diagram:
\begin{equation*}
\xymatrix@1{
&&& \rhand \ar@{<-}[r]^{\phantom{=}=} & \rhand \ar@/^/[dd]^G \ar@/^/@{->}[dll]^\Phi  \\
\diag \ar[r]^\iota_\sim &\diag^S \ar@/_/@{->>}[drr]_\psi \ar@/^/@{->>}[urr]^\phi \ar@{->>}[r]^\pi & \qdiagS
\ar@/^/@{->>}[dr]_{\overline{\psi}}
\ar@/_/@{->>}[ur]^{\overline{\phi}} &&\\
&&& \rstl \ar@{<-}[r]_{\phantom{i}=} & \rstl \ar@/^/[uu]^F \ar@/_/@{->}[ull]_\Psi}
\end{equation*}

At this stage, we do not know whether the functors $\overline{\phi}$ and $\Phi$ (resp.\@ $\overline{\psi}$ and $\Psi$) are
isomorphisms and so inverse each other.

\subsection{Proof of Theorem~\ref{rstl2diag}}\label{proofrstl2diag}
Before proving Theorem~\ref{rstl2diag}, we first establish some lemmas.

\begin{figure}[ht!]\anchor{simulmult}
\begin{center}
      \psfrag{E}[cc][cc]{\small$\Psi_0(P_\ell)$}
      \psfrag{D}[cc][cc]{\small$\Psi_0(P)$}
      \psfrag{i}[Bc][Bc][.5]{$i$}
      \psfrag{m}[Bc][Bc][.5]{$i+1$}
      \psfrag{w}[Bc][Bc][.5]{$i-1$}
      \scaledraw{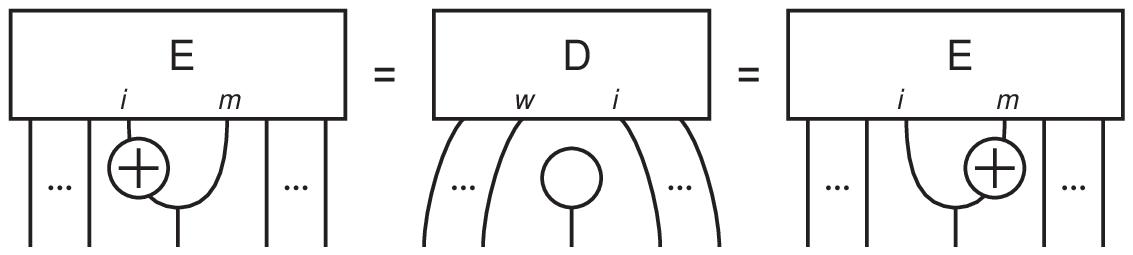}
\end{center}
   \nocolon\caption{}
   \label{simulmult}
\end{figure}

\begin{lem}\label{lemPsi1}
Let $P$ be a ribbon pure $n$-braid and $1 \leq i \leq n+1$. Insert a new component $\ell$ between the $(i-1)$-th and $i$-th
strand of $P$ so that $P\cup \ell$ is a ribbon pure $(n+1)$-braid. Let $P_\ell=\Delta_\ell(P \cup \ell)$ be the ribbon pure
$(n+2)$-braid from $P \cup \ell$ by doubling $\ell$. Then:
\begin{enumerate}
  \renewcommand{\labelenumi}{{\rm (\alph{enumi})}}
  \item The equalities of \figref{simulmult} hold;
  \item $C_i\bigl(\Psi_0(P_\ell)D\bigr)=\Psi_0(P) C_i(D)$  and $C_{i+1}\bigl(D\Psi_0(P_\ell)\bigr)=C_{i+1}(D)\Psi_0(P)$
        for any Hopf diagram $D$ with $(n+2)$ inputs.
\end{enumerate}
\end{lem}
\begin{proof}
Let us prove Part (a) by induction on the length $m$ of $P\cup l$ in the generators $\sigma_{k,l}^{\pm 1}$ and $t_k^{\pm 1}$
of $\RPB_{n+1}$. If $m=0$, then it is an immediate consequence of \eqref{counit} and \eqref{antip3}. Suppose that $m=1$.
Given a ribbon pure braid $Q$, denote by $\Delta_i(Q)$ (resp.\@ $\delta_i(Q)$) the ribbon pure braid obtained from $Q$ by
doubling (resp.\@ deleting) its $i$-th component. We have to verify that the statement is true for
$P=\delta_i(\sigma_{k,l}^{\pm 1})$ and $P_\ell=\Delta_i(\sigma_{k,l}^{\pm 1})$, and for $P=\delta_i(t_k^{\pm 1})$ and
$P_\ell=\Delta_i(t_k^{\pm 1})$. This can be done case by case by using the descriptions of \tabref{tabdesc}. Examples of
such verifications are depicted in \figref{demsimulmultA1}.
\begin{table}[h]\anchor{tabdesc}
   \begin{center}
\begin{tabular}{|c|c|c|}
 \cline{2-3}
 \multicolumn{1}{c|}{} & $\Delta_i(\sigma_{k,l}^{\pm 1})$ & $\delta_i(\sigma_{k,l}^{\pm 1})$ \TabTable \\
 \hline
 $i<k$ & $\sigma_{k+1,l+1}^{\pm 1}$ & $\sigma_{k-1,l-1}^{\pm 1}$ \TabTable \\
 \hline
 $i=k$ & $(\sigma_{i,l+1} \sigma_{i+1,l+1})^{\pm 1}$ & $I_n$ \TabTable \\
 \hline
 $k<i<l$ & $\sigma_{k,l+1}^{\pm 1}$ & $\sigma_{k,l-1}^{\pm 1}$ \TabTable \\
 \hline
 $i=l$ & $(\sigma_{k,i} \sigma_{k,i+1})^{\pm 1}$ & $I_n$ \TabTable \\
 \hline
 $l<i$ & $\sigma_{k,l}^{\pm 1}$ & $\sigma_{k,l}^{\pm 1}$ \TabTable \\
 \hline
\end{tabular}
\;
\begin{tabular}{|c|c|c|}
 \cline{2-3}
 \multicolumn{1}{c|}{} & $\Delta_i(t_k^{\pm 1})$ & $\delta_i(t_k^{\pm 1})$ \TabTable \\
 \hline
 $i<k$ & $t_{k+1}^{\pm 1}$ & $t_{k-1}^{\pm 1} $ \TabTable \\
 \hline
 $i=k$ & $\sigma_{i,i+1}^{\pm 1} t_i^{\pm 1}t_{i+1}^{\pm 1}$ & $I_n$ \TabTable \\
 \hline
 $k<i$ & $t_k^{\pm 1} $ & $t_k^{\pm 1} $  \TabTable \\
 \hline
\end{tabular}
   \end{center}
   \nocolon\caption{}
   \label{tabdesc}
\end{table}

\begin{figure}[ht!]
\begin{center}
      \psfrag{E}[cc][cc]{\small$\Psi_0(\delta_i(\sigma_{i,l}^{-1}))$}
      \psfrag{D}[cc][cc]{\small$\Psi_0(\Delta_i(\sigma_{i,l}^{-1}))$}
      \psfrag{A}[cc][cc]{\small$\Psi_0(\Delta_i(t_i))$}
      \psfrag{B}[cc][cc]{\small$\Psi_0(\delta_i(t_i))$}
      \psfrag{,}[Bl][Bl]{,}
      \psfrag{.}[Bl][Bl]{.}
      \psfrag{i}[Bc][Bc][.5]{$i$}
      \psfrag{m}[Bc][bc][.5]{$i+1$}
      \psfrag{n}[Bc][Bc][.5]{$i-1$}
      \psfrag{w}[Bc][Bc][.5]{$l+1$}
      \scaledraw{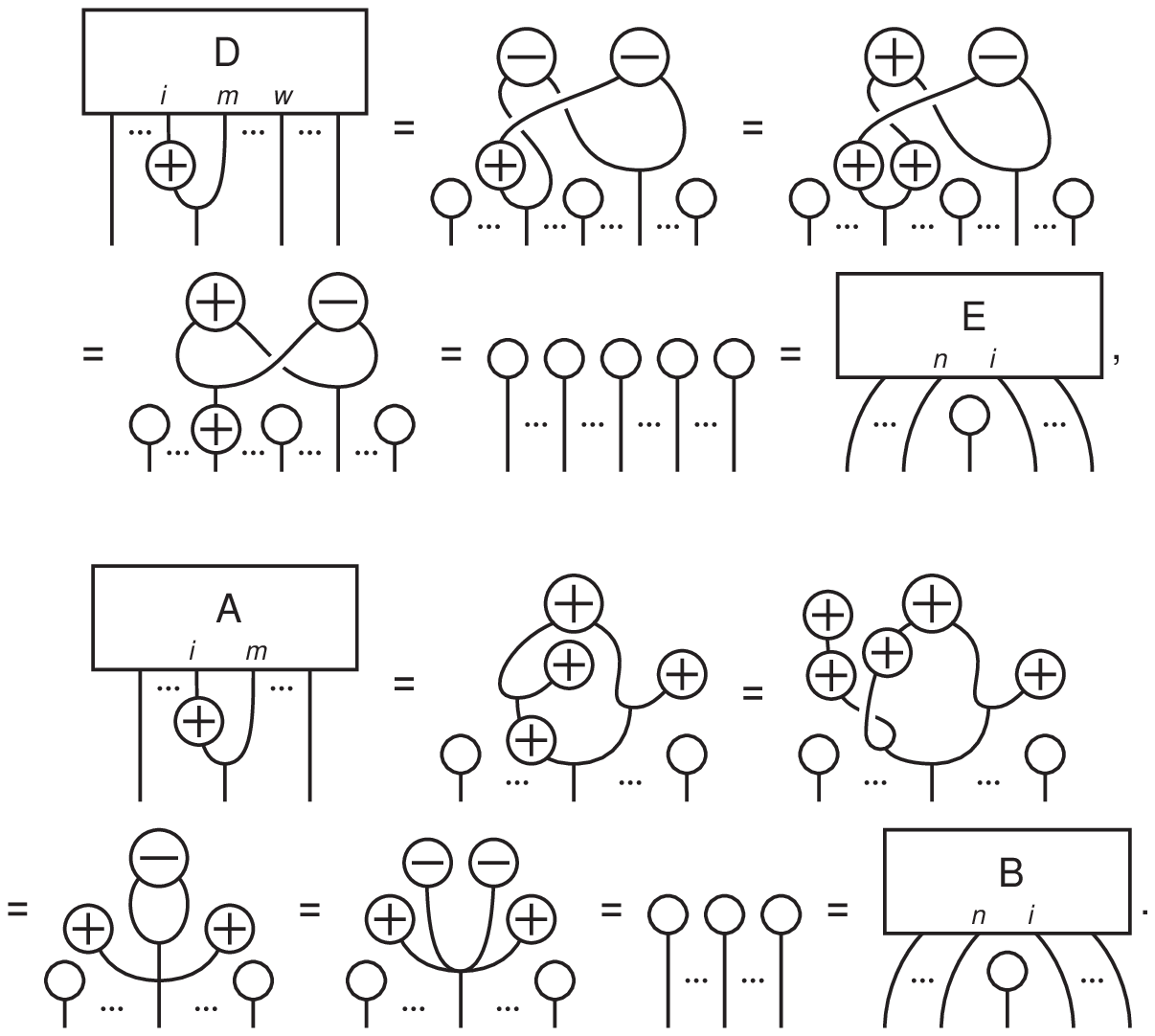}
\end{center}
   \nocolon\caption{}
   \label{demsimulmultA1}
\end{figure}

Let $m\geq 1$ and suppose the statement true for rank $m$. Let $P$ and $\ell$ such that $P \cup \ell=w_1 \dots w_{m+1}$
where the $w_j$ are generators of $\RPB_{n+1}$. Remark that $P_\ell=\Delta_i(w_1) \Delta_i(w_2 \cdots w_{m+1})$ and
$P=\delta_i(w_1) \delta_i(w_2 \cdots w_{m+1})$. By using \eqref{antip2}, \eqref{counit}, and the statement for ranks $1$ and
$m$, we get the equalities depicted in \figref{demsimulmultA2}. The left equality of \figref{simulmult} is then true
for $P$ and $\ell$. The right equality can be verified similarly by remarking that $P_\ell=\Delta_i(w_1\cdots w_m)
\Delta_i(w_{m+1})$. Hence the statement is true for rank $m+1$.
\begin{figure}[ht!]\anchor{demsimulmultA2}
\begin{center}
      \psfrag{E}[cc][cc][.7]{$\Psi_0(\Delta_i(w_1))$}
      \psfrag{D}[cc][cc][.7]{$\Psi_0(\Delta_i(w_2 \cdots w_{m+1}))$}
      \psfrag{F}[cc][cc][.7]{$\Psi_0(\delta_i(w_1))$}
      \psfrag{C}[cc][cc][.7]{$\Psi_0(\delta_i(w_2 \cdots w_{m+1}))$}
      \psfrag{A}[cc][cc]{\small$\Psi_0(P_\ell)$}
      \psfrag{B}[cc][cc]{\small$\Psi_0(P)$}
      \psfrag{i}[Bc][Bc][.5]{$i$}
      \psfrag{m}[Bc][Bc][.5]{$i+1$}
      \psfrag{n}[Bc][Bc][.5]{$i-1$}
      \scaledraw{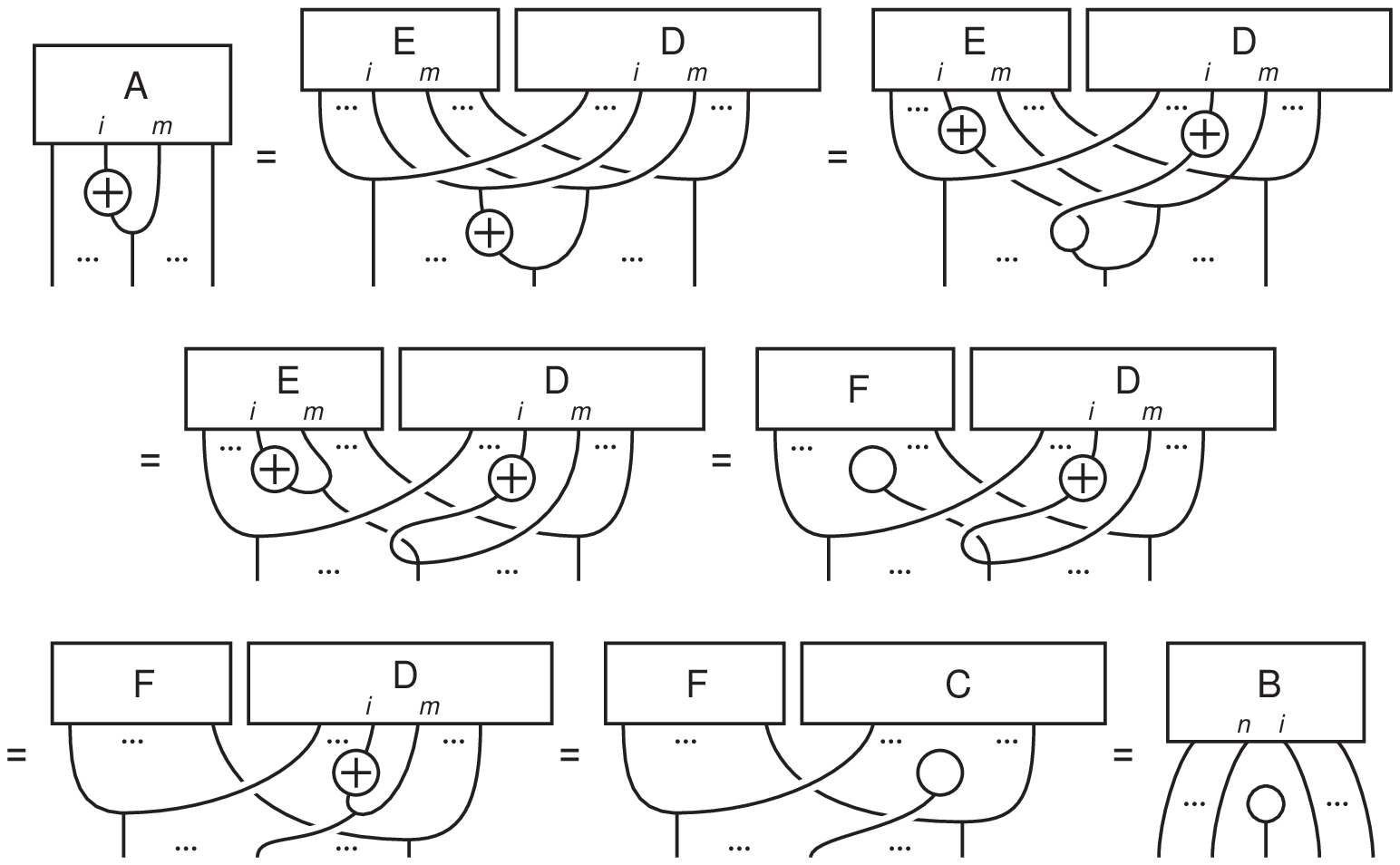}
\end{center}
   \nocolon\caption{}
   \label{demsimulmultA2}
\end{figure}

Let us prove Part (b). Let $D$ be a Hopf $(n+2)$-diagram. By using \eqref{comult}, \eqref{counit}, \eqref{antip2}, and Part
(a) of the lemma, we get the equalities of \figref{demsimulmultB1}, which mean that
$C_i\bigl(\Psi_0(P_\ell)D\bigr)=\Psi_0(P) C_i(D)$. Likewise, one can show that
$C_{i+1}\bigl(D\Psi_0(P_\ell)\bigr)=C_{i+1}(D)\Psi_0(P)$.
\begin{figure}[ht!]\anchor{demsimulmultB1}
\begin{center}
      \psfrag{E}[cc][cc]{\small$\Psi_0(P_\ell)$}
      \psfrag{P}[cc][cc]{\small$\Psi_0(P)$}
      \psfrag{D}[cc][cc]{\small$D$}
      \psfrag{B}[cc][cc]{\small$C_i(D)$}
      \psfrag{i}[cc][cc][.5]{$i$}
      \psfrag{m}[cc][cc][.5]{$i+1$}
      \psfrag{w}[cc][cc][.5]{$i-1$}
      \scaledraw{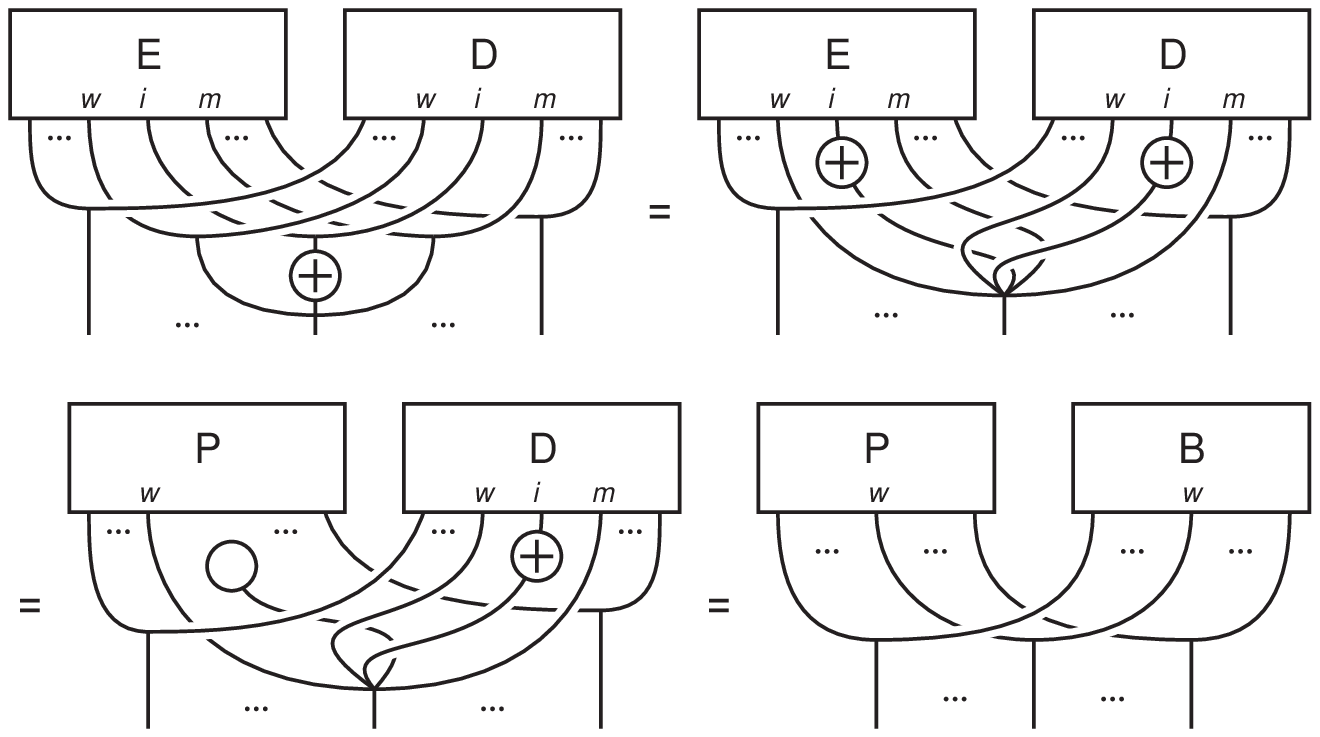}
\end{center}
   \nocolon\caption{}
   \label{demsimulmultB1}
\end{figure}
\end{proof}

\begin{figure}[ht!]\anchor{lemindeppull}
\begin{center}$P=$
      \scaleraisedraw{.45}{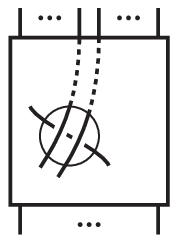} and $P'=$ \scaleraisedraw{.45}{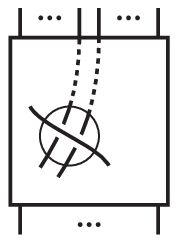}\;
      $\Bigl ($resp.\@ $P=$ \rotatebox[origin=c]{180}{\scaleraisedraw{.45}{psidemo2a}} and
      $P'=$ \rotatebox[origin=c]{180}{\scaleraisedraw{.45}{psidemo2b}} $\Bigr)$
\end{center}
   \nocolon\caption{}
   \label{lemindeppull}
\end{figure}

\begin{lem}\label{lemPsi2}
Let $n\geq 3$ and $1<i<n$ (resp.\@ $0 < i < n-1$). Let $P$ and $P'$ be two pure $n$-braids which have diagrams which differ
only inside a disk. Inside this disk, the $i$-th and $(i+1)$-th strands pass respectively to the front and the back of
another strand. Suppose also that above (resp.\@ below) the disk, the $i$-th and $(i+1)$-th strands run parallel, see
\figref{lemindeppull}. Then:
\begin{equation*}
C_i\Psi_0(P)=C_i\Psi_0(P') \text{ (resp.\@ } C_{i+1}\Psi_0(P)=C_{i+1}\Psi_0(P')\text{).}
\end{equation*}
\end{lem}
\begin{proof}
Let us denote by $k$ the other component. Suppose that above the disk, the $i$\trait th and $(i+1)$-th strands run parallel
(the resp.\@ case can be done similarly). Assume first that $k<i$. We can write $P=A \sigma_{k,i}\sigma_{k,i+1} B$ and $P'=A
B$ where $A$ and $B$ are ribbon pure $n$-braids. Moreover $A=\Delta_\ell(E \cup \ell)$ where $E$ is ribbon pure braid with
$n-2$ strands and $\ell$ is a new component inserted in $E$ in $i$-th position. Now remark that
$\sigma_{k,i}\sigma_{k,i+1}=\Delta_\gamma(I_{n-2} \cup \gamma)$ where $I_{n-2}$ is the trivial braid with $n-2$ strands and
$\gamma$ is new component added to $I_{n-2}$ in $i$-th position so that $I_{n-2} \cup \gamma=\sigma_{k,i}$ in $\RPB_{n-1}$.
Therefore, by using Lemma~\ref{lemPsi1}, we get:
\begin{align*}
C_i\Psi_0(P)&=C_i\Psi_0(A \sigma_{k,i}\sigma_{k,i+1} B)=C_i\Psi_0\bigl(\Delta_\ell(E \cup
\ell)\Delta_\gamma(I_{n-2} \cup \gamma)B\bigr)\\
 & = \Psi_0(E) \, \Psi_0(I_{n-2}) \, C_i\Psi_0(B)= \Psi_0(E) \, C_i\Psi_0(B)\\
 &= C_i\Psi_0\bigl(\Delta_\ell(E \cup\ell)B\bigr)
 = C_i\Psi_0(A B)= C_i\Psi_0(P').
\end{align*}
The case $k>i+1$ is done similarly by writing $P=A \sigma_{i,k}\sigma_{i+1,k} B$ and $P'=AB$.
\end{proof}

Let us now prove Theorem~\ref{rstl2diag}. Let $n$ be a non-negative integer and $T$ be a ribbon $n$\trait string link.
Consider a diagram $\Gamma_T$ of $T$. Applying rules \eqref{renderhanded} changes (in a unique manner) $\Gamma_T$ to a left
handed diagram. Denote by $\alpha_i$ is the number of modifications \eqref{renderhanded} made on the $i$-th component of
$\Gamma_T$. Then pulling maxima to the top and minima to the bottom as explained in Section~\ref{algorstl2diag} leads to a
diagram $P$ of a ribbon pure braid so that $T=(t_1^{\alpha_1} \cdots t_n^{\alpha_n}) c_{j_m} \cdots c_{j_1} \{P\}$ for some
$1 \leq j_1 \leq \dots \leq j_m \leq n$. Pulling extrema in another way may lead to another diagram $P'$ of a ribbon pure
braid so that $T=(t_1^{\alpha_1} \cdots t_n^{\alpha_n})c_{j_m} \cdots c_{j_1} \{P'\}$. Now the ribbon pure braids $\{P\}$ and
$\{P'\}$ are related by moves described in \figref{lemindeppull}. Therefore, by using Lemmas~\ref{contraccomm} and
\ref{lemPsi2}, we obtain: $C_{j_m} \cdots C_{j_1} \Psi_0 (\{P\})=C_{j_m} \cdots C_{j_1} \Psi_0 (\{P'\})$. Hence
\begin{equation*}
\Psi_{\Gamma_T}= (\Omega_1^{\alpha_1} \cdots \Omega_n^{\alpha_n} ) \,  C_{j_m} \cdots C_{j_1} \Psi_0 (\{P\})
\end{equation*}
only depends on the diagram $\Gamma_T$ of $T$.

Let us verify that $\Psi_{\Gamma_T}$ remains unchanged when applying to $\Gamma_T$ a Reidemeister's move, see
\figref{reidemeister}. Invariance under moves of type 2 or type 3 is a consequence of the existence of the functor
$\Psi_0$. Invariance under moves of type 2' is a consequence of Lemmas~\ref{contraccomm} and \ref{lemPsi2}.

\begin{figure}[t]\anchor{figcasD}
   \begin{center}
       \subfigure[]{$\Gamma_T=$ \scaleraisedraw{.45}{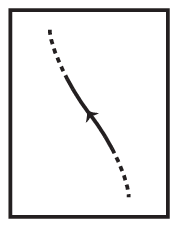} $\;\mapsto \;\Gamma'_T=$ \scaleraisedraw{.45}{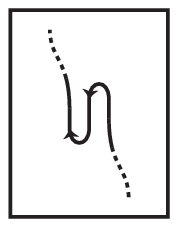}
                    $\;\rightsquigarrow\;$ \scaleraisedraw{.45}{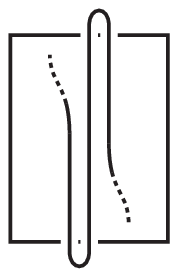}
                    \label{figcasA}}\\
       \subfigure[]{$\Gamma_T=$ \scaleraisedraw{.45}{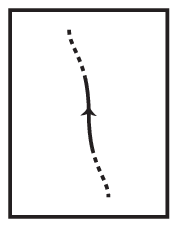} $\;\mapsto \;\Gamma'_T=$ \scaleraisedraw{.45}{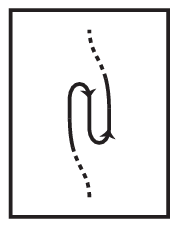}
                    $\;\rightsquigarrow\;$ \scaleraisedraw{.45}{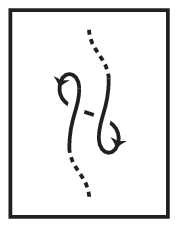} $\;\rightsquigarrow\;$ \scaleraisedraw{.45}{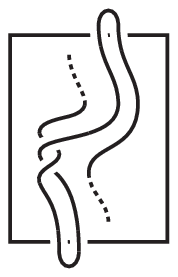}
                    \label{figcasB}}\\
       \subfigure[]{$\Gamma_T=$ \scaleraisedraw{.45}{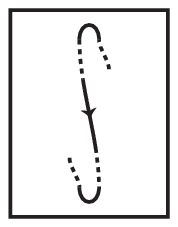} $\;\mapsto \;\Gamma'_T=$ \scaleraisedraw{.45}{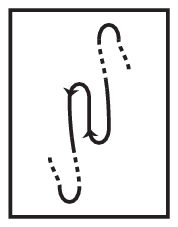}
                    $\;\rightsquigarrow\;$ \scaleraisedraw{.45}{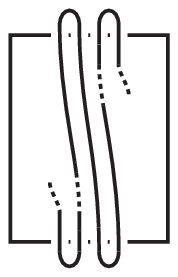}
                    \label{figcasC}} \\
       \subfigure[]{$\Gamma_T=$ \scaleraisedraw{.45}{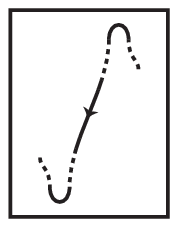} $\;\mapsto \;\Gamma'_T=$ \scaleraisedraw{.45}{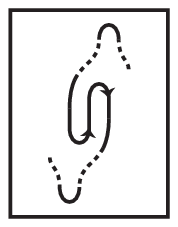}
                    $\;\rightsquigarrow\;$ \scaleraisedraw{.45}{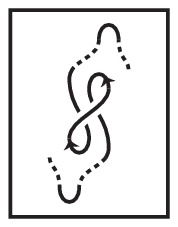} $\;\rightsquigarrow\;$ \scaleraisedraw{.45}{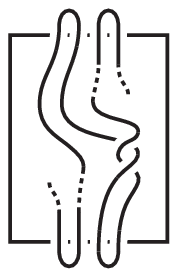}
                    \label{figcasD}}
   \end{center}
   \nocolon\caption{}
   \label{figcas}
\end{figure}

\begin{figure}[t]\anchor{figcasparticD}
   \begin{center}
\begin{minipage}[c]{.92\textwidth}
\psfrag{i}[cc][cc][.8]{$i$} \psfrag{k}[cc][cc][.8]{$k$} \psfrag{U}[cc][cc]{\small$U$} \psfrag{V}[cc][cc]{\small$V$} $\Gamma_T=
\scaleraisedraw{.45}{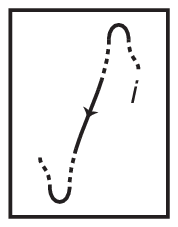} \;\rightsquigarrow \; \scaleraisedraw{.45}{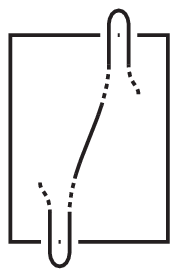} \;\rightsquigarrow \;c_{j_m} \cdots
c_{j_1} \Bigl ( \,\scaleraisedraw{.54}{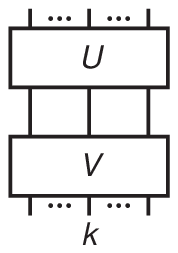} \,\Bigl
)$\\
\hspace*{1.47cm}\rotatebox{270}{$\mapsto$}\\
\psfrag{U}[cl][cl]{\small$U_\ell$} \psfrag{V}[cl][cl]{\small$V_{\ell'}$} $\Gamma'_T= \scaleraisedraw{.45}{casD2}
\;\rightsquigarrow\;\scaleraisedraw{.45}{casD4}  \;\rightsquigarrow\; c_{j_m} \cdots c_{j_r} c_{i+1} c_{j_{r-1}}\cdots
c_{j_1} \Bigl ( \,\scaleraisedraw{.54}{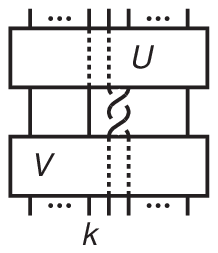} \,\Bigl )$
\end{minipage}
   \end{center}
   \nocolon\caption{}
   \label{figcasparticD}
\end{figure}

\begin{figure}[t]\anchor{grossefig}
\begin{center}
      \scalebox{.95}{\psfrag{E}[cc][cc]{\small$\Psi_0(U_\ell)$}
      \psfrag{D}[cc][cc]{\small$\Psi_0(V_{\ell'})$}
      \psfrag{P}[cc][cc]{\small$\Psi_0(U)$}
      \psfrag{Q}[cc][cc]{\small$\Psi_0(V)$}
      \psfrag{R}[cc][cc]{\small$\Psi_0(U)\Psi_0(V)$}
      \psfrag{k}[Bc][Bc][.5]{$k-1$}
      \psfrag{w}[Bc][Bc][.5]{$k$}
      \psfrag{i}[Bc][Bc][.5]{$k+1$}
      \psfrag{m}[Bc][Bc][.5]{$k+2$}
      \psfrag{t}[Bc][Bc][.5]{$k+3$}
      \scaledraw{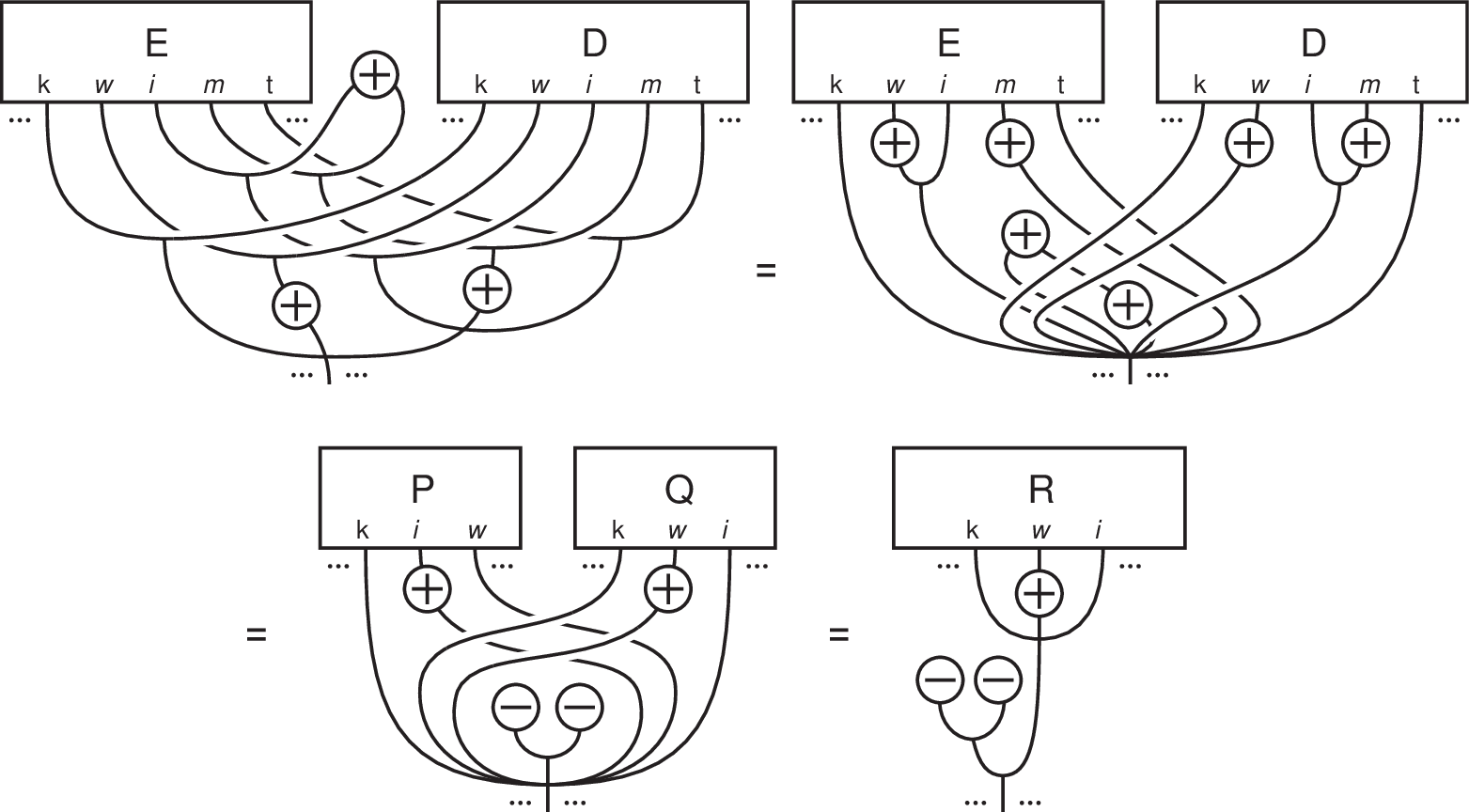}}
\end{center}
   \nocolon\caption{}
   \label{grossefig}
\end{figure}

Suppose that a Reidemeister move of Type 0 is applied to $\Gamma_T$, and denote $\Gamma'_T$ the diagram so obtained. There
are four cases to consider, depending on the orientation of the considered strand and the direction (left or right handed) of
the move. These cases, together with the way we apply the algorithm, are depicted in \figref{figcas}.
Let us for example verify invariance in the case depicted in \figref{figcasD}. Recall that applying the algorithm to
$\Gamma_T$ gives rise to a diagram $P$ of a pure braid such that $\{\Gamma_T\} =(t_1^{\alpha_1} \cdots t_n^{\alpha_n})  \,
c_{j_m} \cdots  c_{j_1}  \{ P \}$. Let $i$ be the number of the component of $\Gamma_T$ on which the move is performed.
Since the orientation of the strand is downwards and $T$ is a (canonically oriented) string link, we know that there exist a
maximum just before and a minimum just after the place where the move is performed, when going through the $i$\trait th
component from bottom to top (see the left picture of \figref{figcasD}). Let $c_{j_r}$ be the contraction corresponding
to this pair of extrema when applying the algorithm to $\Gamma_T$ (we have $j_r=i+1$). Denote by $k$ the number of the
component of $P$ where the move is performed. Up to using invariance under Reidemeister's moves of type~2 and type 3, we can
write $\{P\}=UV$, with $U$ and $V$ pure braids, so that the move is performed  on the $k$-th component between $U$ and $V$.
Insert a new component $\ell$ in $U$ in $k$-th position and a new component $\ell'$ in $V$ in $(k+1)$-th position. Set
$U_\ell=\Delta_\ell(U \cup \ell)$ and $V_{\ell'}=\Delta_{\ell'}(V \cup \ell')$. As depicted in Figures~\fref{figcasD}
and~\fref{figcasparticD}, we can apply the algorithm to $\Gamma'_T$ in such a way that:
\begin{equation*}
\{\Gamma'_T\} =(t_1^{\alpha_1} \cdots t_i^{\alpha_i+2}\cdots t_n^{\alpha_n})  c_{j_m} \cdots c_{j_r} c_{i+1}
c_{j_{r-1}}\cdots c_{j_1} (U_\ell \sigma_{k+1,k+2} V_{\ell'}).
\end{equation*}
Now, by using \eqref{antip2}, \eqref{antip4} and \eqref{relsuptwist1}, we have, for any Hopf diagram $D$,
\begin{equation}\label{ciomegp}
C_j (\Omega_p D )= \left \{ \begin{array}{ll} \Omega_p C_j (D) &\text{for $p<j$,}\\ \Omega_{j-1} C_j (D) & \text{for
$p=j$,}\\ \Omega_{p-2} C_j (D) & \text{for $p>j$.} \end{array} \right.
\end{equation}
Moreover we have the equalities of \figref{grossefig} where are used in particular \eqref{antip5}, \eqref{relsuptwist1},
\eqref{relsuptwist4} and Lemma~\ref{lemPsi1}(a). Then we get that:
\begin{align*}
C_kC_k \Psi_0(U_\ell \sigma_{k+1,k+2} V_{\ell'}) &= C_kC_{k+2} \bigl ( \Psi_0(U_\ell) \Sigma_{k+1,k+2}
\Psi_0(V_{\ell'})\bigr )\\
&=\Omega_{k-1}^{-2} C_k \bigl (\Psi_0(U) \Psi_0(V) \bigr )\\
&=\Omega_{k-1}^{-2} C_k\Psi_0(\{P\}).
\end{align*}
Therefore we can conclude that:
\begin{align*}
\psi_{\Gamma'_T} &= (\Omega_1^{\alpha_1} \cdots \Omega_i^{\alpha_i+2}\cdots \Omega_n^{\alpha_n} ) \, C_{j_m}
\cdots C_{j_r} C_{i+1} C_{j_{r-1}}\cdots C_{j_1}\Psi_0(U_\ell \sigma_{k+1,k+2} V_{\ell'}) \\
&= (\Omega_1^{\alpha_1} \cdots \Omega_i^{\alpha_i+2}\cdots \Omega_n^{\alpha_n} ) \, C_{j_m}
\cdots \widehat{C_{j_r}}\cdots C_{j_1} C_k C_k \Psi_0(U_\ell \sigma_{k+1,k+2} V_{\ell'}) \\
&=  (\Omega_1^{\alpha_1} \cdots \Omega_i^{\alpha_i+2}\cdots \Omega_n^{\alpha_n}) \, C_{j_m}
\cdots \widehat{C_{j_r}} \cdots C_{j_1} \bigl (\Omega_{k-1}^{-2} C_k\Psi_0(\{P\}) \bigr )\\
&= (\Omega_1^{\alpha_1} \cdots \Omega_i^{\alpha_i+2}\cdots \Omega_n^{\alpha_n})  \Omega_i^{-2} \, C_{j_m}
\cdots \widehat{C_{j_r}} \cdots C_{j_1} C_k \Psi_0(\{P\}) \quad \text{by \eqref{ciomegp}}\\
&= (\Omega_1^{\alpha_1} \cdots \Omega_n^{\alpha_n})\, C_{j_m} \cdots C_{j_1} \Psi_0(\{P\}) =\psi_{\Gamma_T}.
\end{align*}
Invariance in the cases depicted in Figures~\fref{figcasA}, \fref{figcasB}, and \ref{figcasC} can be checked as above.

By using similar techniques, one can show the invariance of $\Psi_{\Gamma_T}$ under Reidemeister's moves of type 1'. Hence
we conclude that $\Psi_{\Gamma_T}$ remains unchanged when applying to $\Gamma_T$ a Reidemeister's move, and so
$\Psi(T)=\Psi_{\Gamma_T}$ is well-defined.

We get directly from its construction that the functor $\Psi$ is monoidal and satisfies Condition (a). Let us check that it
satisfies Condition (b). Let $1<i < n$ and $T$ be a ribbon $n$-string link with trivial $i$-th component. By applying the
algorithm we have $T=(t_1^{\alpha_1} \cdots t_n^{\alpha_n}) c_{j_m} \cdots c_{j_1} (P)$ for some $1 \leq j_1 \leq \dots \leq
j_m \leq n$ and some pure braid~$P$. Since the $i$-th component of $T$ is trivial, we can apply the algorithm to a diagram
of $c_i(T)$ in such a way that:
\begin{equation*}
c_i(T)=(t_1^{\alpha_1} \cdots t_{i-2}^{\alpha_{i-2}} t_{i-1}^{\alpha_{i-1}+\alpha_i+\alpha_{i+1}} t_i^{\alpha_{i+2}}\cdots
t_{n-2}^{\alpha_n}) c_{j_m-2}\cdots c_{j_{k+1}-2} c_i c_{j_k} \cdots c_{j_1} (P),
\end{equation*}
where $k$ is such that $j_{k+1}> i+1 \geq j_k$. Hence:
\begin{align*}
C_i(\Psi &(T)) = C_i \bigl (\Omega_1^{\alpha_1}\cdots \Omega_n^{\alpha_n} C_{j_m} \cdots C_{j_1}
\Psi_0(P) \bigr )\\
&=(\Omega_1^{\alpha_1} \cdots \Omega_{i-1}^{\alpha_{i-1}+\alpha_i+\alpha_{i+1}}
\cdots \Omega_{n-2}^{\alpha_n})  C_i C_{j_m} \cdots C_{j_1} \Psi_0(P) \quad \text{by \eqref{ciomegp}}\\
&=(\Omega_1^{\alpha_1} \cdots \Omega_{i-1}^{\alpha_{i-1}+\alpha_i+\alpha_{i+1}}
\cdots \Omega_{n-2}^{\alpha_n}) C_{j_m-2}\cdots C_{j_{k+1}-2} C_i C_{j_k} \cdots C_{j_1} \Psi_0(P)\\
& = \Psi(c_i(T)).
\end{align*}
Uniqueness of a functor $\qdiagS \to \rstl$ satisfying (a) and (b) comes from the fact that every ribbon string link can be
realized as a sequence of contractions of a ribbon pure braid such that each of contraction is performed on a trivial
component (see the algorithm described in Section~\ref{algorstl2diag}).

Finally, let $T$ be a ribbon string link. We can always write $T=c_{j_m} \cdots c_{j_1}(P)$ for some ribbon pure braid $P$.
Then we have:
\begin{align*}
\overline{\psi} \Psi(T) &= \overline{\psi} \Psi(c_{j_m} \cdots c_{j_1}(P)) \\
  & = \overline{\psi} \bigl (  C_{j_m} \cdots C_{j_1}\Psi_0(P) \bigr ) \quad \text{by Conditions (a) and (b)} \\
  & = c_{j_m} \cdots c_{j_1} \bigl (\overline{\psi} \Psi_0(P) \bigr ) \quad \text{by Lemma~\ref{lemcontr}} \\
  & = c_{j_m} \cdots c_{j_1} (P)=T \quad \text{by Lemma~\ref{proppure2diag}.}
\end{align*}
Hence $\overline{\psi}\circ\Psi=1_{\rstl}$. This completes the proof of Theorem~\ref{rstl2diag}.

\section{Quantum invariants via Hopf diagrams and Kirby elements}\label{sect-qinv}

A general method is given in \cite{Vir} for defining quantum invariants of 3\trait manifolds starting from a ribbon category
(or a ribbon Hopf algebra). In this section, we explain the role played by Hopf diagrams in this theory. Note that this was
the initial motivation of this work.

\subsection{Dinatural transformations and coends}\label{dinatcoend} We give here definitions adapted to our purposes.
For more general situations, we refer to \cite{ML1}.

Let $\cc$ be a category with left duals. By a \emph{dinatural transformation of $\cc$}, we shall mean a pair $(Z,d)$
consisting in an object $Z$ of $\cc$ and a family $d$, indexed by $\Ob(\cc)$, of morphisms $d_X\co \dual{X} \otimes X \to Z$
in $\cc$ satisfying $d_Y(\id_{\dual{Y}} \otimes f)=d_X(\dual{f} \otimes \id_X)$ for any morphism $f\co  X \to Y$.

By a \emph{coend of $\cc$}, we shall mean a dinatural transformation $(A,i)$ which is universal in the sense that, if
$(Z,d)$ is any dinatural transformation, then there exists a unique morphism $r\co A \to Z$ such that $d_X=r \circ i_X$ for
all object $X$ in $\cc$, see \figref{figcoend}. Note that a coend, if it exists, is unique (up to unique isomorphism).
\begin{figure}[ht!]\anchor{figcoend}
$$
\xymatrix{
\dual{Y} \otimes X \ar[d]_{\dual{f} \otimes 1}\ar[r]^{1 \otimes f}  & \dual{Y} \otimes Y \ar[d]_{d_Y} \ar@/^/[rdd]^{i_Y}& \\
\dual{X} \otimes X \ar[r]^{d_X} \ar@/_/[rrd]_{i_X} & Z \ar@{<.}[rd]|-{\exists ! r} &\\
&& A}
$$
     \nocolon\caption{}
     \label{figcoend}
\end{figure}

\begin{rems}
(1)\qua In general, the coend of $\cc$ always exists in a completion of $\cc$, namely the category of $\mathrm{Ind}$-objects of
$\cc$ (see \cite{Lyu1}). However, for simplicity,  we will restrict to the case where the coend exists in~$\cc$.

(2)\qua If $\cc$ is the category of representations of a finite-dimensional Hopf algebra or is a premodular category, then the
coend exists in $\cc$, see~\cite{Vir}.
\end{rems}
Assume that $(A,i)$ is the coend of $\cc$. Let the morphisms $\ev_X\co \dual{X} \otimes X \to \un$ and $\coev_X\co \un \to X
\otimes \dual{X}$ be the evaluation and coevaluation associated to the left dual $\dual{X}$ of an object $X$. The following
dinatural transformations:
\begin{equation*}
(i_X \otimes i_X)(1 \otimes \coev_X \otimes 1)\co  \dual{X} \otimes X \to A \otimes A \quad \text{and} \quad \ev_X\co
\dual{X} \otimes X \to \un
\end{equation*}
factorize respectively to morphisms $\Delta_A \co  A \to A \otimes A$ (the \emph{coproduct}) and $\varepsilon_A\co A \to \un$
(the \emph{counit}). This makes $A$ a coalgebra in the category~$\cc$.

\subsection{The coend of a ribbon category}\label{ribboncoend} Let $\cc$ be a ribbon category (see \cite{Tur2}) and assume
that the coend $(A,i)$ of $\cc$ exists. We denote the braiding of $\cc$ by $c_{X,Y}\co X\otimes Y \to Y \otimes X$, and the
twist of $\cc$ by $\theta_X\co X \to X$.

The coalgebra structure $(\Delta_A,\varepsilon_A)$ on the object $A$ (see Section~\ref{dinatcoend}) extends to a structure
of a Hopf algebra, see \cite{Lyu1}. This means that there exist morphisms $\mu_A \co  A \otimes A \to A$ (the product),
$\eta_A\co  \un \to A$ (the unit), and $S_A \co  A \to A$ (the \emph{antipode}). They satisfy the same axioms as those of a
Hopf algebra except the usual flip is replaced by the braiding $c_{A,A}\co  A \otimes A \to A \otimes A$. Namely, the unit
morphism is $\eta_A= i_\un\co \un = \dual{\un} \otimes \un \to A$. By using the universal property of a coend\footnote{For
the product $\mu_A$, use it twice.}, the product $\mu_A$ and the antipode $S_A$ are defined as follows:
\begin{align*}
& \mu_A(i_X \otimes i_Y)=i_{Y \otimes X} (1_{\dual{X}} \otimes c_{X,\dual{Y} \otimes Y})\co  \dual{X} \otimes X \otimes
\dual{Y}
\otimes Y \to A,\\
& S_A i_X=(\ev_X \otimes i_{\dual{X}})(\id_{\dual{X}} \otimes c_{\leftidx{^{\vee\vee}}{\! X}{},X} \otimes
\id_{\dual{X}})(\coev_{\dual{X}} \otimes c_{\dual{X},X})\co \dual{X} \otimes X \to A.
\end{align*}
Note that $S_A$ is invertible, with inverse $S_A^{-1}\co A \to A$ defined via:
\begin{equation*}
S_A^{-1} i_X=(\ev_X \otimes i_{\dual{X}})(\id_{\dual{X}} \otimes c^{-1}_{\leftidx{^{\vee\vee}}{\! X}{},X} \otimes
\id_{\dual{X}})(\coev_{\dual{X}} \otimes c^{-1}_{\dual{X},X}) \co \dual{X} \otimes X \to A.
\end{equation*}
Let us define the morphisms $\omega_A\co  A \otimes A \to \un$ and $\theta^\pm_A\co A \to \un$ as follows:
\begin{align*}
&\theta_A^\pm i_X=\ev_X(\id_{\dual{X}} \otimes \theta_X^{\pm 1})\co \dual{X} \otimes X  \to \un,\\
&\omega_A (i_X \otimes i_Y)=\omega_{X,Y}\co \dual{X} \otimes X \otimes \dual{Y} \otimes Y \to \un,
\end{align*}
where $\omega_{X,Y}=(\ev_X \otimes \ev_Y)(\id_{\dual{X}} \otimes c_{\dual{Y},X}c_{X,\dual{Y}} \otimes \id_{\dual{Y}})$. It
can be shown that~$\omega_A$ is a Hopf pairing, see \cite{Vir}. Finally, we set:
\begin{equation*}
\omega_A^+=\omega_A(S_A^{-1} \otimes \id_A) \quad \text{and} \quad \omega_A^-=\omega_A.
\end{equation*}

\subsection{Hopf diagrams and factorization}\label{Hopffacto} Let $\cc$ be a ribbon category. Assume that the coend $(A,i)$
of $\cc$ exists. Consider the Hopf algebra structure of $A$ and the morphisms $\omega_A^\pm$ and $\theta_A^\pm$ as in
Section~\ref{ribboncoend}.

Let $T$ be a ribbon $n$-handle. Recall that for any objects $X_1, \dots, X_n$ of $\cc$, the handle $T$ defines a morphism
$T_{X_1, \dots, X_n}\co \dual{X_1} \otimes X_1 \otimes \cdots \otimes \dual{X_n} \otimes X_n \to \un$ by canonically
orienting $T$ (see Section~\ref{sectrhand}) and decorating the $k^\mathrm{th}$ component of $T$ by~$X_k$. Hence it can be
factorized thought the coend to a morphism $T_\cc\co A^{\otimes n} \to \un$ so that:
\begin{equation}\label{phiTL}
T_{X_1, \dots, X_n}=T_\cc \circ (i_{X_1} \otimes \cdots \otimes i_{X_n})
\end{equation}
for all objects $X_1, \dots, X_n$ of $\cc$.

Let $\qddS$ be the braided category defined in Section~\ref{reladdhopfdiag}.

\begin{thm}\label{Evaldiag}
\begin{enumerate}
\renewcommand{\labelenumi}{{\rm (\alph{enumi})}}
  \item There exists a unique braided functor $e_\cc \co\qddS \to \cc$ sending the object $*$ to $A$ and
the morphisms $\Delta$, $\varepsilon$, $S^{\pm 1}$,  $\omega_\pm$, $\theta_\pm$ to the morphims $\Delta_A$, $\varepsilon_A$,
$S_A^{\pm 1}$, $\omega_A^\pm$, $\theta^\pm_A$ respectively.
  \item The functor $e_\cc$ induces a monoidal functor $E_\cc\co \qdiagS \to \conv_\cc(A,\un)$.
  \item For any ribbon $n$-handle
$T$, the factorization morphism $T_\cc\co A^{\otimes n} \to \un$ defined by $T$ is given by $T_\cc=E_\cc(D)$ where $D$ is
any Hopf diagram such that $\overline{\phi}(D)=T$. In particular, $T_\cc=E_\cc\circ \Phi(T)$.
\end{enumerate}
\end{thm}
\begin{rem}
We expect to generalize this result to categories~$\cc$ which are not ribbon but only turban, see \cite{Brug3}.
\end{rem}

\begin{proof}
Let us prove Part (a). Uniqueness of the functor $e_\cc$ is clear since we impose the image of all its generators. Its
existence comes from the fact that the relations imposed on the generators $\Delta$, $\varepsilon$, $S^{\pm 1}$,
$\omega_\pm$, $\theta_\pm$ and $\tau_{*,*}^{\pm 1}$ in $\qddS$ are still true in $\cc$ when replacing them  by $\Delta_A$,
$\varepsilon_A$, $S_A^{\pm 1}$, $\omega_A^\pm$, $\theta^\pm_A$ and $c_{A,A}^{\pm 1}$ respectively (since $\omega_A$ is a Hopf
pairing and $\theta_A$ is defined using the twist of $\cc$).

Part (b) is a direct consequence of $\qdiagS=\conv_{\qddS}(*,\un)$.

Finally, let us prove Part (c). By the uniqueness of the factorization morphism via a coend, we have to show that if $D$ is
a Hopf diagram with $n$ inputs, then
\begin{equation}\label{demphiTL}
\overline{\phi}(D)_{X_1, \dots, X_n}=E_\cc(D) \circ (i_{X_1} \otimes \cdots \otimes i_{X_n})
\end{equation}
for any objects $X_1, \dots, X_n$ of $\cc$. Clearly, if \eqref{demphiTL} is true for two Hopf diagrams $D$ and $D'$, then it
is also true for the Hopf diagram $D\otimes D'$.

Moreover, by definition of $\phi$ (see \figref{Hopf2Hand}) and of $\varepsilon_A$, $\omega_A^\pm$ and $\theta^\pm_A$,
\eqref{demphiTL} is true for the trivial Hopf diagram $\varepsilon \otimes \cdots \otimes \varepsilon$ and for the Hopf
diagrams $\Sigma_{i,j}^{\pm 1}$ and $\Omega_k^{\pm 1}$ depicted in \eqref{gammaomega}. Now suppose that \eqref{demphiTL} is
true for some Hopf diagram $D$. Then,  by definition of $\phi$ (see \figref{Hopf2Hand}) and of $\Delta_A$, $S_A^{\pm 1}$
and $c_{A,A}^{\pm 1}$, \eqref{demphiTL} remains true for the following diagrams:
\begin{equation*}
      \psfrag{D}[cc][cc]{\small$D$}\scaleraisedraw{.47}{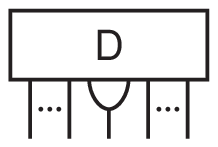}\; , \quad
      \psfrag{D}[cc][cc]{\small$D$}\scaleraisedraw{.47}{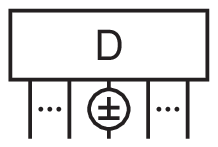}\; , \quad
      \psfrag{D}[cc][cc]{\small$D$}\scaleraisedraw{.47}{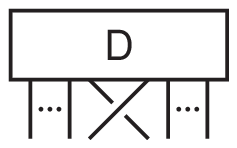}\; , \quad \text{and} \quad
      \psfrag{D}[cc][cc]{\small$D$}\scaleraisedraw{.47}{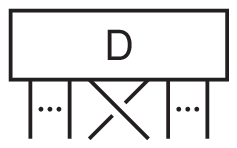}\; .
\end{equation*}
Hence we can deduce that \eqref{demphiTL} is always true.
\end{proof}

\subsection{Kirby elements} Let $\cc$ and $A$ be as in Section~\ref{Hopffacto}. Let $L$ be a framed link in~$S^3$
with $n$ components. Choose an orientation for $L$. There always exists a ribbon $n$\trait handle $T$ (not necessarily
unique) such that $L$ is isotopic $T \circ (\cup \otimes \cdots \otimes \cup)$, where $\cup$ denotes the cup with
counterclockwise orientation and $T$ is canonically oriented. For $\alpha \in \Hom_\cc(\un,A)$, set:
\begin{equation}
\tau_\cc(L;\alpha)=T_\cc \circ \alpha^{\otimes n} \in \End_\cc(\un),
\end{equation}
where $T_\cc\co A^{\otimes n} \to \un$ is defined as in \eqref{phiTL}. Following \cite{Vir}, by a \emph{Kirby element}
of~$\cc$ we shall mean a morphism $\alpha \in \Hom_\cc(\un,A)$ such that, for any framed link~$L$, $\tau_\cc(L;\alpha)$ is
well-defined and invariant under isotopies and 2-handle slides of~$L$.

In general, determining the set of the morphisms $T_\cc$ when $T$ runs over ribbon handles is quite difficult. Nevertheless,
by Theorem~\ref{Evaldiag}, the $T_\cc$'s belong to the set of morphisms given by evaluations of Hopf diagrams. Hence a
criterion for a morphism to be a Kirby element:

\begin{cor}\label{corKirby}
Let $\alpha\co  \un \to A$ in $\cc$. Suppose that the two following conditions are satisfied:
\begin{enumerate}
  \renewcommand{\labelenumi}{{\rm (\alph{enumi})}}
  \item for any integer $n\geq 1$ and any Hopf diagram $D$ with $n$ entries, we have:
\begin{equation*}
        \psfrag{a}[cc][cc]{\small$\alpha$}
       \psfrag{D}[cc][cc]{\small$E_\cc(D)$}
       \psfrag{S}[cc][cc]{\small$S_A$}
       \scaleraisedraw{.5}{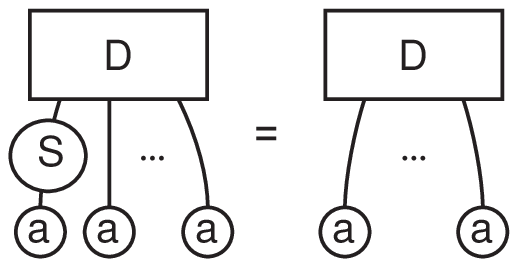}\quad;
\end{equation*}
  \item for any integer $n\geq 2$ and any Hopf diagram $D$ with $n$ entries, we have:
\begin{equation*}
       \psfrag{a}[cc][cc]{\small$\alpha$}
       \psfrag{D}[cc][cc]{\small$E_\cc(D)$}
       \scaleraisedraw{.5}{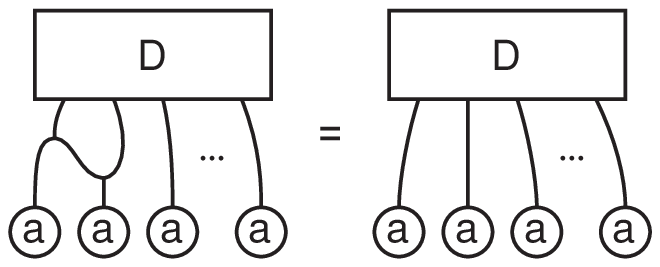} ,
\end{equation*}
where \scaleraisedraw{.11}{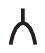} and \scaleraisedraw{.11}{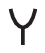} denote $\mu_A\co A \otimes A \to A$ and $\Delta_A\co
A \to A \otimes A$ respectively.
\end{enumerate}
Then $\alpha$ is a Kirby element of $\cc$.
\end{cor}

\begin{rem}
One recovers the fact (see \cite[Theorem~2.5]{Vir}) that a morphism $\alpha\co  \un \to A$ is a Kirby element if it
satisfies:
\begin{equation*}
\left \{
\begin{array}{l}
S_A \alpha=\alpha,\\
(\mu_A \otimes \id_A)(\id_A \otimes \Delta_A)(\alpha \otimes \alpha)=\alpha \otimes \alpha,
\end{array}
\right.
\end{equation*}
or, in case $\cc$ is linear, if the morphisms
\begin{equation*}
\left \{
\begin{array}{l}
S_A \alpha-\alpha,\\
(\mu_A \otimes \id_A)(\id_A \otimes \Delta_A)(\alpha \otimes \alpha)-\alpha \otimes \alpha,
\end{array}
\right.
\end{equation*}
are negligible.
\end{rem}
\begin{proof}
We adapt the proof of \cite[Theorem~2.5]{Vir} to this more general situation. Let $L=L_1 \cup \dots \cup L_n$ be a framed
link. Choose an orientation of $L$ and a ribbon $n$-handle $T$ such that $L$ is isotopic $T \circ (\cup \otimes \cdots
\otimes \cup)$.

Suppose firstly that we reverse the orientation of $L_i$. Let us denote by $L'$ the oriented framed link obtained from $L$
by this orientation change. As in the proof of \cite[Theorem~2.5]{Vir}, we can choose a ribbon $n$-handle $T'$ such $L'$ is
isotopic to $T' \circ (\cup \otimes \cdots \otimes \cup)$ and $T'_\cc=T_\cc \circ (\id_{A^{\otimes (i-1)}} \otimes S_A
\otimes \id_{A^{\otimes (n-i)}})$. Now let $D$ be a Hopf diagram such that $\overline{\phi}(D)=T$. Then $T_\cc=E_\cc(D)$ and
so
\begin{align*}
\tau_\cc(L';\alpha) & = E_\cc(D) (\id_{A^{\otimes (i-1)}} \otimes S_A \otimes \id_{A^{\otimes (n-i)}}) \alpha^{\otimes n}\\
   & = E_\cc(D) (c_{A,A^{\otimes (i-1)}} \otimes \id_{A^{(n-i)}})(S_A \alpha \otimes \alpha^{\otimes n-1})\\
   & = E_\cc (D') (S_A \alpha \otimes \alpha^{\otimes n-1})\\
   & = E_\cc (D')\alpha^{\otimes n} \quad \text{by Condition (a)}\\
   & = E_\cc(D)(c_{A,A^{\otimes (i-1)}} \otimes \id_{A^{(n-i)}})\alpha^{\otimes n}\\
   & = E_\cc(D)\alpha^{\otimes n}=\tau_\cc(L;\alpha),
\end{align*}
where
\begin{equation*}
    D'=  \psfraga <2pt,0pt> {D}{\small$D$}\psfrag{i}[B][B][.7]{\small$i$}\scaleraisedraw{.47}{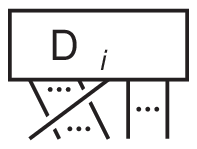} \; .
\end{equation*}
Hence $\tau_\cc(L;\alpha)$ does not depend on the choice of an orientation for $L$.

Suppose now that the component $L_2$ slides over the component $L_1$. Let $L'_1$ be a parallel copy of $L_1$ and set $L'=L_1
\cup (L'_1 \# L_2) \cup L_3 \cup \dots \cup L_n$. As in the proof of \cite[Theorem~2.5]{Vir}, we can choose a ribbon
$n$-handle $T'$ such that $L'$ is isotopic to $T' \circ (\cup \otimes \cdots \otimes \cup)$ and $T'_\cc=T_\cc \circ \bigl
((\mu_A \otimes \id_A)(\id_A \otimes \Delta_A) \otimes \id_{A^{\otimes (n-2)}} \bigr )$. Hence, choosing a Hopf diagram $D$
such that $\overline{\phi}(D)=T$, we get that $T_\cc=E_\cc(D)$ and so
\begin{align*}
\tau_\cc(L';\alpha)& = E_\cc(D) \bigl ((\mu_A \otimes \id_A)(\id_A \otimes \Delta_A) \otimes \id_{A^{\otimes (n-2)}} \bigr )
\alpha^{\otimes n}\\
 &= E_\cc(D) \alpha^{\otimes n} \quad \text{by Condition (b)}\\
 &=\tau_\cc(L;\alpha).
\end{align*}
Likewise, using the braiding,  we can show that $\tau_\cc(L;\alpha)$ is invariant under the other handle slides. Hence
$\alpha$ is a Kirby element of $\cc$.
\end{proof}

\subsection{Quantum invariants of $3$-manifolds} Recall (see \cite{Li}) that every closed, connected and oriented
3-dimensional manifold can be obtained from $S^3$ by surgery along a framed link $L \subset S^3$.

For any framed link $L$ in $S^3$, we will denote by $S^3_L$ the 3-manifold obtained from $S^3$ by surgery along $L$, by
$n_L$ the number of components of $L$, and by $b_-(L)$ the number of negative eigenvalues of the linking matrix of $L$.

A Kirby element $\alpha$ of $\cc$ is said to be \emph{normalizable} if $\theta^+_A \alpha$ and $\theta^-_A \alpha$ are
invertible in the semigroup $\End_\cc(\un)$. To each normalizable Kirby element $\alpha$ is associated an invariant
$\tau_\cc(M;\alpha)$ of (closed, connected, and oriented) 3-manifolds $M$ with values in  $\End_\cc(\un)$, see
\cite[Proposition~2.3]{Vir}. From Theorem~\ref{Evaldiag}, we immediately deduce that:

\begin{cor}\label{dercoro}
For a normalizable Kirby element $\alpha$ and any framed link $L$, we have
\begin{equation*}
\tau_\cc(S^3_L;\alpha)=(\theta^+_A \alpha)^{b_-(L)-n_L}\,  (\theta^-_A \alpha)^{-b_-(L)} \; E_\cc(D) \alpha^{\otimes n},
\end{equation*}
where $D$ is any Hopf diagram such that $L$ is isotopic to $\phi(D) \circ (\cup \otimes \cdots \otimes \cup)$.
\end{cor}

\begin{rems}
1) Corollary~\ref{dercoro} gives an intrinsic description, in terms of Hopf algebraic structures, of quantum invariants of
3-manifolds.

2) In Corollary~\ref{dercoro}, we can in particular take $D=\Psi(T)$ where $\Psi$ is as in Theorem~\ref{rstl2diag} and $T$
is a ribbon sting link such that $L$ is isotopic to closure of $T$. Recall also that the constructive proof of
Theorem~\ref{rstl2diag} provides us with an algorithm for computing $D=\Psi(T)$ starting from a diagram of $T$.
\end{rems}

Let us conclude by giving an example. Recall that the Poincar\'e  sphere $\mathbb{P}$ can be presented as $\mathbb{P}\simeq
S^3_K$ where $K$ is the right-handed trefoil with framing $+1$. Now this trefoil can be represented as $K=\phi(D) \circ
\cup$ where $D=(\omega_+ \Delta \otimes \theta_-)\Delta$. Indeed, we have:
\begin{equation*}
  K=\scaleraisedraw{.45}{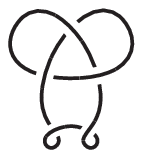} \quad \text{and} \quad D=\scaleraisedraw{.403}{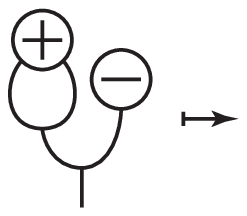}  \;\; \phi(D)= \,
      \psfrag{-}[cc][cc]{\scalebox{1.25}{$\sim$}}\scaleraisedraw{.403}{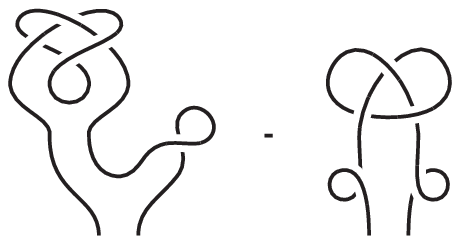}.
\end{equation*}
Hence, if $\alpha\co \un \to A$ is a normalizable Kirby element of $\cc$, then:
\begin{equation*}
 \tau_\cc(\mathbb{P};\alpha)=\;
 \psfrag{o}[c][c]{\small$\theta^+_A$}
 \psfrag{a}[c][c]{\small$\alpha$}
 \psfrag{1}[Bl][Bl]{$-1$}
 \scaleraisedraw{.45}{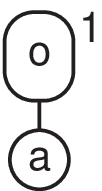} \quad \;\,
 \psfrag{D}[c][c]{\small$\Delta_A$}
 \psfrag{o}[c][c]{\small$\theta^-_A$}
 \psfrag{w}[c][c]{\small$\omega^+_A$}
 \psfrag{a}[c][c]{\small$\alpha$}
 \scaleraisedraw{.45}{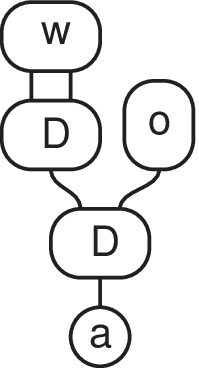}.
\end{equation*}

\Addresses\recd

\end{document}